\documentclass{article}

\usepackage{arxiv/arxiv}

\usepackage[utf8]{inputenc} 
\usepackage[T1]{fontenc}    
\usepackage{hyperref}       
\usepackage{url}            
\usepackage{booktabs}       
\usepackage{amsfonts}       
\usepackage{nicefrac}       
\usepackage{microtype}      
\usepackage{graphicx}
\usepackage[numbers]{natbib}
\usepackage{doi}
\usepackage{bbm}
\usepackage{amsmath}
\usepackage{subfig}
\usepackage{orcidlink}
\usepackage{enumitem}
\newlength{\jeroenlen}

\newcommand{\zh}{\mathbf{z}_h}
\newcommand{\yh}{\mathbf{y}_h}
\newcommand{\xh}{\mathbf{x}_h}
\newtheorem{remark}{Remark}
\newcommand{\tildeyh}{\tilde{\mathbf{y}}_h}
\newcommand{\yhrec}{\mathbf{y}_{h,\mathrm{rec}}}
\newcommand{\yn}{\mathbf{y}_N}
\newcommand{\tildeyn}{\tilde{\mathbf{y}}_N}
\newcommand{\uh}{\mathbf{u}_h}
\newcommand{\tildeuh}{\tilde{\mathbf{u}}_h}
\newcommand{\uhrec}{\mathbf{u}_{h,\mathrm{rec}}}
\newcommand{\un}{\mathbf{u}_N} 
\newcommand{\tildeun}{\tilde{\mathbf{u}}_N}
\newcommand{\mus}{\boldsymbol{\mu}_s}

\newcommand{\R}{\mathbb{R}}

\newcommand{\citeinline}[2]{\cite{#1}}
\newcommand{\citesinline}[4]{\cite{#1,#2}}

\newcommand{\citesssinline}[8]{\cite{#1,#2,#3,#4}}
\newcommand{\myspace}{\]\[}

\title{Real-time optimal control of \\ high-dimensional parametrized systems \\ by deep learning-based reduced order models}

\author{\href{https://orcid.org/0000-0002-7111-5070}{\includegraphics[scale=0.06]{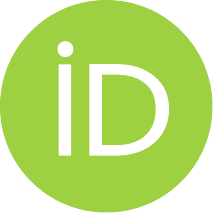}\hspace{1mm}Matteo Tomasetto} \\
	Department of Mechanical Engineering\\
	Politecnico di Milano, Milano, Italy \\
        \texttt{matteo.tomasetto@polimi.it}
	\And 
 \href{https://orcid.org/0000-0001-8277-2802}{\includegraphics[scale=0.06]{orcid.pdf}\hspace{1mm}Andrea Manzoni} \\
        MOX - Department of Mathematics \\
        Politecnico di Milano, Milano, Italy \\
        \texttt{andrea1.manzoni@polimi.it}
	\AND
 \href{https://orcid.org/0000-0002-0476-4118}{\includegraphics[scale=0.06]{orcid.pdf}\hspace{1mm}Francesco Braghin } \\
	Department of Mechanical Engineering\\
	Politecnico di Milano, Milano, Italy \\
       \texttt{francesco.braghin@polimi.it}
        }

\date{}

\renewcommand{\headeright}{}
\renewcommand{\undertitle}{}
\renewcommand{\shorttitle}{Real-time optimal control of high-dimensional parametrized systems by DL-ROMs}

\hypersetup{
pdftitle={Real-time optimal control of high-dimensional parametrized
systems by DL-ROMs},
pdfauthor={Matteo Tomasetto, Andrea Manzoni, Francesco Braghin}
}

\begin{document}
\maketitle

\begin{abstract}
Steering a system towards a desired target in a very short amount of time is a challenging task from a computational standpoint. Indeed, the intrinsically iterative nature of optimal control problems requires multiple simulations of the state of the physical system to be controlled. Moreover, the control action needs to be updated whenever the underlying scenario undergoes variations, as it often happens in applications. Full-order models based on, e.g., the Finite Element Method, do not  meet these requirements due to the computational burden they usually entail. On the other hand, conventional reduced order modeling techniques such as the Reduced Basis method, despite their rigorous construction, are intrusive, rely on a linear superimposition of modes, and lack of efficiency when addressing nonlinear time-dependent dynamics. In this work, we propose a non-intrusive Deep Learning-based Reduced Order Modeling (DL-ROM) technique for the rapid control of systems described in terms of parametrized PDEs in multiple scenarios. In particular, optimal full-order snapshots are generated and properly reduced by either Proper Orthogonal Decomposition or deep autoencoders (or a combination thereof) while feedforward neural networks are exploited to learn the map from scenario parameters to reduced optimal solutions. Nonlinear dimensionality reduction therefore allows us to consider state variables and control actions that are both low-dimensional and distributed. After {(i)} data generation, {(ii)} dimensionality reduction, and {(iii)} neural networks training in the offline phase, optimal control strategies can be rapidly retrieved in an online phase for any scenario of interest. The computational speedup and the extremely high accuracy obtained with the proposed approach are finally assessed on different PDE-constrained optimization problems, ranging from the minimization of energy dissipation in incompressible flows modelled through Navier-Stokes equations to the thermal active cooling in heat transfer.
\end{abstract}

\keywords{optimal control \and PDE-constrained optimization \and parametrized systems \and reduced order modeling \and deep learning \and fluid flows \and heat transfer}

\section{Introduction}\label{sec:introduction}


Controlling dynamical systems is one of the most challenging and widespread tasks in Applied Sciences and Engineering. For instance, reducing energy consumption, ensuring system stability, or attaining a desired configuration, are just three examples of problems that can be cast in the framework of Optimal Control Problems (OCPs). In order to steer the dynamics of a physical system towards a desired target, a suitable control action has to be properly selected and tuned. The optimal control strategy is typically found as a minimum point of a \textit{loss} or \textit{cost functional}. Minimization has to be constrained by the governing equation driving the considered dynamics in order to narrow the space of possible solutions on those that are physically admissible \cite{Manzoni2021}. The governing equations may be, in general, algebraic systems, ordinary differential equations or, as considered in the following, Partial Differential Equations (PDEs). Specifically, we take into account PDEs parametrized by a vector of input (or scenario) parameters $\mus$. Indeed, our main goal is to conceive an efficient numerical strategy to determine optimal controls in multiple scenarios of interest as fast as possible. OCPs governed by PDEs are traditionally solved by exploiting high-fidelity full-order models (FOMs) based on discretization techniques such as, e.g., the Finite Element Method (FEM). Despite their ubiquitous use -- from fluid flows, backs to the 90s \cite{gunzburger1991analysis,hinze2001second,gunzburger2002perspectives} to more recent applications in robotic swarms \citep{Sinigaglia2022, Sinigaglia2022a, Sinigaglia2022b, Sinigaglia2022c} and acoustic and thermal cloaking \citep{Cominelli2022, Chen2021, Sinigaglia2022-ROM, Sinigaglia_2024} -- relying on FOMs becomes unfeasible in the case of multiple scenarios, since an iterative optimization procedure and thus several PDE solves would be required for each new scenario.

Reduced Order Modeling (ROM) techniques \citep{Hesthaven2015, Quarteroni2015, Brunton2019, Manzoni2019}, such as the Reduced Basis (RB) method relying on Proper Orthogonal Decomposition (POD), have been largely exploited to reduce the dimensionality of optimality equations and thus enable faster resolutions \citep{Kunisch1999, Kunisch2008, negri2013reduced, Benner2014, Amsallem2015, Sinigaglia2022-ROM}. In the context of Model Predictive Control (MPC),  POD and system identification techniques, such as Dynamic Mode Decomposition (DMD) \citep{Schmid2010}, Sparse Identification of Nonlinear Dynamics (SINDy) \citep{Brunton2016, Brunton2016a} or Eigensystem Realization Algorithm (ERA), have been employed to replace the physical constraint by a low-rank surrogate model \cite{Kaiser2018, Hickner2023, Proctor2016}. Similar strategies are exploited to describe the dynamics of usually low-dimensional observables measured by sensors, that is, an approximation of the Koopman operator \cite{Peitz2019, Klus2020, Korda2018}. However, the linear and intrusive nature of the RB method does not make it suitable, or computationally affordable, whenever {\em (i)} the finite element matrices are not accessible, {\em (ii)} if the physical laws are nonlinear, or {\em (iii)} if the problems are not affine with respect to the scenario parameters.

In addition to ROM strategies, neural networks have been more recently considered in the literature to speed up the resolution of OCPs. For instance, in the MPC framework, the differential equations appearing as constraints are often approximated through neural network-based surrogate models \cite{Draeger1995, Chen2018, Bieker2019, Peitz2023, Antonelo2024}, such as Recurrent Neural Networks, Physics-Informed Neural Networks or Input Convex Neural Networks, which replicate the dynamics of the available data, achieving a remarkable computational speedup. However, since no reduction techniques have been systematically applied in these cases, only low-dimensional states and controls can be taken into account to deal with feasible input and output neural networks dimensions. \\

To overcome the limitations of the previous ROM frameworks and thus allowing for greater flexibility and speedup, in this work we exploit non-intrusive ROM strategies to develop a framework able to solve any high-dimensional parametrized OCP in real-time for multiple scenarios. Differently from the RB method, non-intrusive techniques compress the available data through POD \cite{Hesthaven2018}, deep convolutional autoencoders \cite{Fresca2021, Franco2023}, or a combination thereof \cite{Fresca2022}, while the map going from scenarios to the reduced optimal solution can be modeled through a deep feedforward neural network. A combination of linear dimensionality reduction and neural networks has also been considered by \citeinline{Luo2023}{Luo et al} in the context of robust control, while convolutional autoencoders are taken into account, for instance, by \citeinline{Ishize2023}{Ishize et al} to reduce the dimensionality of uncontrolled state snapshots. In this work, instead, we properly reduce also high-dimensional control variables and we focus on controlled snapshots in order to avoid costly optimization procedures when computing optimal actions related to new scenarios.

OCPs have been increasingly addressed with Physics-Informed Neural Networks (PINNs) \citep{Raissi2019} modeling the space-time dependence of state variables. In addition to the misfit between model outputs and experimental data, PINN also encode the residuals of the physical equations in the loss function to be minimized as regularizing terms  to achieve more consistent results. For example, \citesinline{Wang2021}{ Hwang2022}{Wang et al}{Hwang et al} solve OCPs approximating the control-to-state map through Physics-Informed Deep Operator Networks (PI-DeepONets \cite{Wang2021a, Lu2021}). Moreover, all-at-once strategies based on PINNs, where optimal control, state, and adjoint variables are estimated with neural networks properly trained, generally speaking, minimizing the residuals of first-order optimality conditions, have been also investigated \cite{Mowlavi2023, Barry-Straume2022, Demo2023, Yin2024}. Notably, akin to our approach, \citesinline{Demo2023}{Yin2024}{Demo et al}{Yin et al} consider parametrized PDEs to quickly retrieve parameter-specific optimal control strategies: however, in contrast to our approach, the lack of dimensionality reduction techniques and a data-driven approach result in expensive training phases.

Another interesting and promising method to solve control problems is represented by Reinforcement Learning \citep{Sutton2018} and, more recently, Deep Reinforcement Learning (DRL). This framework aims to learn an optimal policy, that is, an optimal state-to-control map, in order to control a physical system from sensor-based measurements or simulated snapshots. This is achieved thanks to an iterative interaction between the agent -- which has to decide the best possible control action to apply -- and the environment, which is the dynamical system under investigation, possibly described by differential equations. At every episode, the agent applies a time-dependent action on the environment and receives the corresponding state and reward values: thanks to this information, it is possible to improve the policy towards the optimum trying to maximize the reward, i.e., to minimize the cost. This strategy has been exploited in several applications such as, e.g., flow control problems \citep{Ren2021a, Verma2018, Rabault2019} or metamaterial design and cloaking \cite{Shah2021, Mirzakhanloo2020}. DRL algorithms usually deal only with low-dimensional variables due to sample inefficiency and computationally demanding training phases. To restrict the number of agent-environment interactions, \citeinline{Zolman2024}{Zolman et al} exploit SINDy to build efficient and interpretable surrogate models for the environment, the reward and the policy, however only addressing low-dimensional dynamical systems and without considering systems described in terms of PDEs. \citeinline{Ma2018}{Ma et al}, instead, handle high-dimensional observations through autoencoders in order to jointly reduce the problem dimensionality and learning an optimal policy more quickly in the context of fluid jets control for ball games. In this work, we consider similar reduction strategies to deal with high-dimensional state and control variables. Moreover, we focus on parametrized PDEs to easily handle state variability, as recently proposed by \citeinline{Botteghi2024}{Botteghi and Fasel} in the context of DRL, and we exploit an offline-online decomposition as done by \citeinline{Sanchez2018}{Sánchez-Sánchez and Izzo}. However, we do not consider feedback signals; rather, we only require to select the new scenario parameters to infer the corresponding optimal pair, with no need of state measurements online. For a complete overview on DRL see, e.g., \cite{Garnier2021, Vignon2023a, Rabault2020} while, for a general introduction on hybrid approaches that combine machine learning techniques into physical problems, see, e.g., \cite{Brunton2020,San2021,Vinuesa2023,Ren2020}. 

The paper is organized as follows. Section~\ref{sec:OCP} reviews parametrized optimal control problems and intrusive dimensionality reduction techniques, such as the RB method. Section~\ref{sec:OCP-DL-ROM} proposes a nonlinear and non-intrusive reduced order modeling strategy based on deep learning to solve optimal control problems of parametrized PDEs in real-time. Section~\ref{sec:test} shows the performances of the proposed approach on three different applications ranging from flow control to active thermal cooling. Section~\ref{sec:conclusions} discusses some future development ideas about possible extensions of the proposed tool.

\section{Optimal control of parametrized partial differential equations}\label{sec:OCP}

This section briefly introduces the formulation of parametrized OCPs, the high-fidelity full-order model based on the finite element method, and standard dimensionality reduction techniques taken into account to speed up the resolution of OCPs.

\subsection{Mathematical formulation of parametrized optimal control problems}

Optimal control problems aim at finding the best control action capable of steering the dynamics of a (state) system as close as possible to a target configuration. This task is usually achieved by minimizing a cost functional subject to the constraint expressed by the state equation -- this latter consisting of a differential problem -- as follows:
\begin{equation}
\label{eq:OCP}
    \mbox{Given} \ \mus, \quad \mbox{find} \quad  J_h(\xh; \mus) \to \min \ \ \ \text{ s.t. } \ \ \ {\bf G}_h(\xh; \mus) = {\bf 0}, \ \xh \in \mathcal{X}_{ad} 
\end{equation}
where $\xh = (\yh, \uh) \in \R^{N_h^y} \times \R^{N_h^u}$ are the optimization (state and control) variables, possibly subject to additional constraints (here $\mathcal{X}_{ad} = \mathcal{Y}_{ad} \times \mathcal{U}_{ad}$). For the sake of simplicity, we directly take into account the discrete state and control variables, obtained by discretizing the state problem by means of a high-fidelity full-order model (FOM) built through, e.g.,  the finite element method (FEM). In practice, the infinite-dimensional state $y$ and control $u$ are approximated by $y_h \in \mathcal{Y}_h$ and $u_h \in \mathcal{U}_h$, where $\mathcal{Y}_h$ and $\mathcal{U}_h$ are finite-dimensional spaces spanned by a basis of $N_h^y$ and $N_h^u$ elements, respectively -- here $h>0$ denotes the discretization parameter related to the mesh size, so that the smaller $h$, the larger $N_h^y$ and $N_h^u$ and the more accurate and computational expensive the FEM approximation. In this way, we can retrieve the algebraic formulation of the differential equations in terms of $\yh \in \R^{N_h^y}$ and $\uh \in \R^{N_h^u}$, that are the basis expansion coefficients of $y_h$ and $u_h$. Notice that, whenever a nodal basis is chosen, $\yh$ and $\uh$ corresponds to the state and control values at the mesh nodes. For a complete presentation of FEM and optimal control problems, see, e.g., \citesinline{Quarteroni2017}{Manzoni2021}{Quarteroni}{Manzoni et al}. In particular:

\begin{enumerate}
\item 
\textbf{State equation.} The governing equation ${\bf G}_h(\xh; \mus) = {\bf 0}$, with ${\bf G}_h \in \R^{N_h^y}$,   encodes the physical law describing the state dynamics in a region of interest $\Omega$ and, possibly, over a time interval $(0,T)$ for a given final time $T > 0$. In this work, we focus on (nonlinear and/or time-dependent) PDEs parametrized by a vector of $p$ scenario parameters $\mus$ belonging to a parameter space $\mathcal{P} \subset \R^p$. For example, $\mus$ may {\em (i)} represent physical or material properties appearing in the equation, {\em (ii)} account for the geometrical variability of the region where the phenomenon takes place, or {\em (iii)} appear as a result of the parametrization of a given term in the PDE. Note that, in case of time-dependent problems, the time variable $t$ is regarded as an additional scenario parameter. The considered PDEs are thus affected by both the control action $\uh$ and the vector of parameters $\mus$. The former is the quantity to tune in order to steer the dynamics towards the target, the latter represents the scenario variability that we want to address. Indeed, we aim to control the dynamics for multiple scenarios as fast as possible. The governing equation also includes suitable boundary and initial conditions on, respectively, $\partial \Omega \times (0,T]$ and $\Omega \times \{t = 0\}$, in order to deal with well-posed problems and to ensure the existence and uniqueness of an optimal solution that depends continuously on the problem data. 

\item 
\textbf{Cost functional.} $J_h(\xh; \mus): \R^{N_h^y} \times \R^{N^u_h} \to \R$ is a function encoding the objective  to be achieved. The minimum point of $J_h$ corresponds to the optimal state and control variables, also referred to as optimal pair. To guarantee the OCP well-posedness and avoid control strategies that are unfeasible from an energetic standpoint, additional regularization terms are typically added in the cost functional like, e.g., the $L^2$-norm of the control and its gradient \cite{Manzoni2021}. Note that, since the information on the state of the system is often incomplete, a partial observation $\zh = C_h\yh \in R^{N^z_h}$, with $C_h: \R^{N_h^y} \to \R^{N^z_h}$, can replace the state in the definition of $J_h$.
\end{enumerate}


\subsection{Parametrized optimal control problems: from full-order models to reduced order models}
\label{subsec:RB}

Solving PDE-constrained optimization problems is a challenging task due to the interplay between the optimization workflow and the approximation of the underlying PDE. In the unconstrained case (when no constraints are set on $\xh$, that is, $\mathcal{X}_{ad}$ is the entire space) a set of Karush-Kuhn-Tucker (KKT) optimality conditions can be obtained through the Lagrange multipliers' method, and reads as follows:
\begin{equation}
{\bf F}_h(\xh, {\bf p}_h; \mus) = \begin{bmatrix}
{\bf g}_h(\xh; \mus) + \mathbb{C}_h^{\top}(\xh; \mus) {\bf p}_h  \smallskip \\
{\bf G}_h(\xh; \mus)
\end{bmatrix}
 = {\bf 0}
\label{eq:KKT_1}
\end{equation}
where ${\bf p}_h \in \mathbb{R}^{N_h^p}$ is the discrete adjoint vector,  
\[
\mathbb{C}_h^{\top}(\xh; \mus) = [
\partial_y {\bf G}_{h} (\xh;\mus) \qquad 
\partial_u {\bf G}_{h} (\xh;\mus)
] \in \R^{N_h^p \times (N_h^y + N_h^u)}
\]
is the state operator Jacobian, whereas
\[
{\bf g}_h  = 
[
\nabla_y J_h (\xh;\mus) \qquad \nabla_u J_h (\xh;\mus)
]^{\top} \in \R^{N_h^y + N_h^u}
\]
is the gradient of $J_h$ with respect to $\xh = (\yh, \uh)$. Assuming that the same space is used to discretize both state and adjoint variables, we have $N_h^p = N_h^y$. To better highlight the structure of the KKT optimality conditions and the role of the adjoint equation, we can equivalently rewrite system  \eqref{eq:KKT_1} as follows:
\begin{equation}
     \label{eq:nonlin_numer_optim}
\left\{
\begin{array}{ll}
\nabla_y J_h ((\yh, \uh);\mus) 
 + (\partial_y {\bf G}_{h} ((\yh, \uh);\mus))^{\top }\mathbf{p}_h= {\bf 0}& \text{ \ \ (adjoint equation)}  \smallskip \\
 \nabla_u J_h ((\yh, \uh);\mus) + (\partial_u {\bf G}_{h} ((\yh, \uh);\mus))^{\top }{\bf p}_h  ={\bf 0}   & \text{ \ \ (optimality condition)} \smallskip \\
 {\bf G}_h\left( (\yh, \uh); \mus \right) ={\bf 0} &  \text{ \ \ (state equation)}. \smallskip%
\end{array}%
\right.
\end{equation}

The optimality system \eqref{eq:nonlin_numer_optim} is usually a coupled, nonlinear (and possibly time-dependent) system that can in principle be solved through the Newton method -- this corresponds to the so-called {\em direct method} in optimal control. Alternatively, the {\em indirect method} takes into account a two-step process where {\em (i)} the solution of the state and the adjoint problems and {\em (ii)} the minimization of the cost functional by means of a descent method  are considered iteratively until convergence. 

\begin{remark}
In the simplest case where a quadratic cost functional 
\[
J_h(\yh, \uh; \mus) = \frac{1}{2} (\yh - {\bf y}_d)^{\top} \mathbb{M}_h (\yh - {\bf y}_d) + \frac{\beta}{2} \uh^\top \mathbb{M}_h \uh
\]
and a linear, stationary state problem 
\[
{\bf G}_h(\yh, \uh; \mus) = \mathbb{A}_h(\mus)\yh - {\bf f}_h(\mus) - \mathbb{B}_h(\mus)\uh
\]
are considered, with $\mathbf{y}_d$ standing for the target state, \eqref{eq:nonlin_numer_optim} becomes a linear system featuring a saddle-point structure, of the following form:
\begin{equation}
\label{eq:KKT_block13}
\underbrace{\begin{bmatrix}
 \   \mathbb{M}_h  \ & \ 0 \  & \ \mathbb{A}^{\top}_h \   \\
 \ 0 \ & \ \beta \mathbb{M}_h \   & \ -\mathbb{M}_h \  \\
 \  \mathbb{A}_h  \ & \   -\mathbb{M}_h \ & \ 0 \ \\
\end{bmatrix}}_{K_h(\mus)}
\underbrace{\begin{bmatrix}
{\bf y}_h \\
{\bf u}_h \\
{\bf p}_h \\
\end{bmatrix}}_{{\bf q}_h(\mus)} =
\underbrace{\begin{bmatrix}
\mathbb{M}_h {\bf y}_{d} \\
{\bf 0} \\
{\bf f}_h \\
\end{bmatrix}}_{{\bf s}_h(\mus)}.  
\end{equation}
where $\mathbb{A}_h$ denotes the stiffness matrix, $\mathbb{M}_h$ the mass matrix, $\mathbb{B}_h$ the control-to-state matrix, $\mathbf{f}_h$ the source term.
\end{remark}

\begin{figure}
    \centering
    \includegraphics[scale = 0.5]{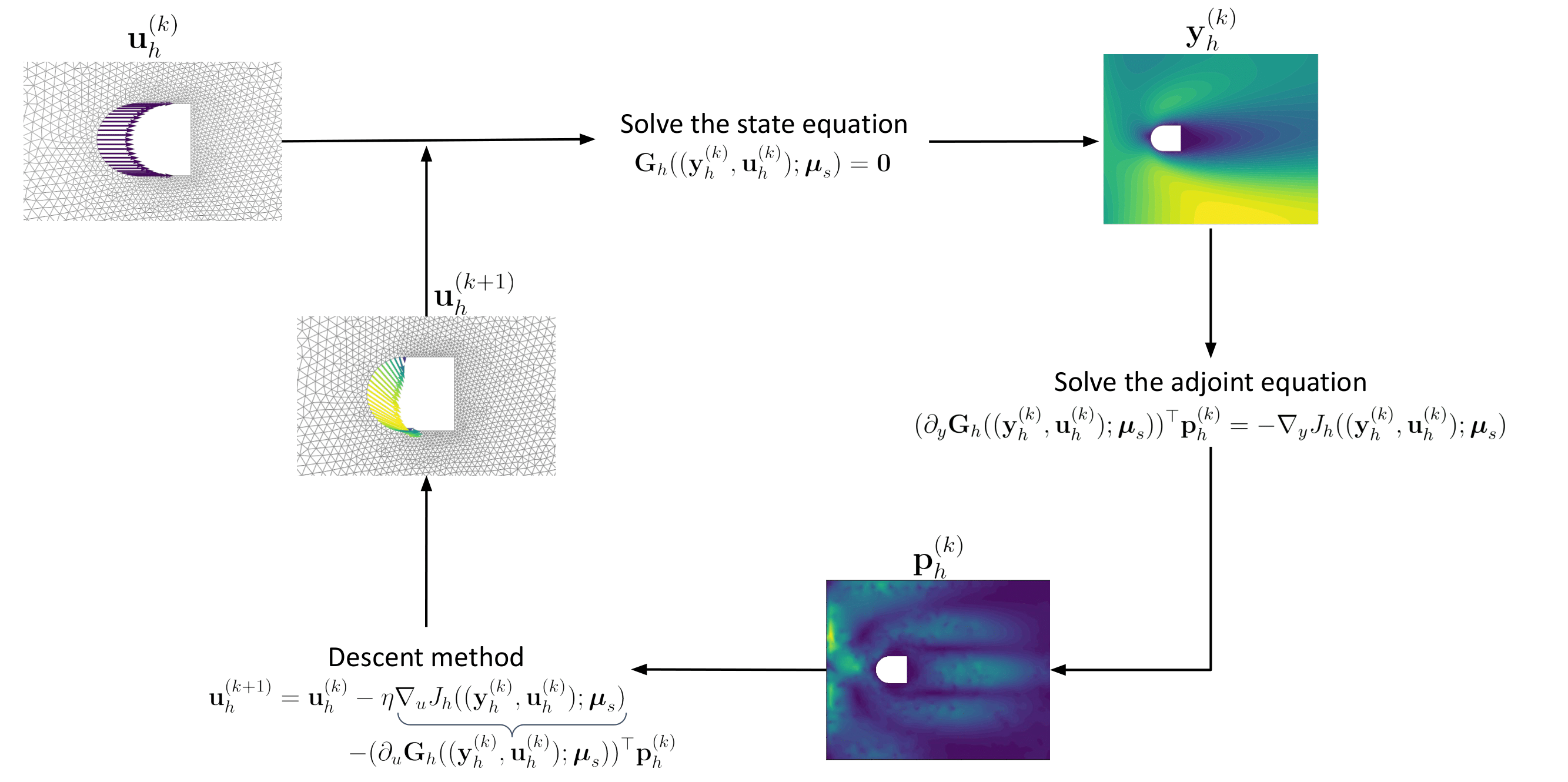}
    \caption{\textit{Test 1.1}. Indirect approach for a PDE-constrained optimization problem in the case of steady fluid flow control. Here $\uh^{(k)}$ is the suboptimal boundary control at the $k$-th iteration of the optimization loop; $\yh^{(k)}$ denotes the corresponding flow velocity computed by solving the steady Navier-Stokes equations; $\mathbf{p}_h^{(k)}$ is the corresponding adjoint variable, resulting from the adjoint equation; $\uh^{(k+1)}$ is the boundary control updated through the steepest descent method with step size $\eta$.}
    \label{fig:OCP_loop}
\end{figure}

Figure~\ref{fig:OCP_loop} shows the rationale behind indirect approaches, focusing on the control problem we will tackle in Section~\ref{subsec:steadyNS} for a steady flow past an obstacle. Starting from an initial control strategy, here expressed as a boundary datum on the obstacle, the flow field is computed together with the adjoint variables in order to compute the gradient $\nabla_{u} J_h(\yh(\uh), \uh)$ and then update the control values toward the minimum of the loss function in case, e.g., a steepest descent method is employed. These steps have to be repeated iteratively until a suitable convergence criterion is met. 

The loop sketched in Figure~\ref{fig:OCP_loop} thus requires a PDE solve at every optimization step in order to compute the state $\yh$ starting from the current control $\uh$. When considering the algebraic formulation obtained through the FEM approximation, this is equivalent to solve a (possibly, nonlinear and/or time-dependent) system of $N_h^y$ equations. Similarly, the adjoint problem requires the solution of a linear (possibly time-dependent) system of $N_h^p$ equations. Especially when taking into account challenging small-scales dynamics, $N_h^y$ (and $N_h^p$) can be rather large, ultimately entailing a huge computational burden. Even more importantly, whenever optimal control strategies must be determined for different scenarios of interest, the entire optimization procedure has to be repeated almost from scratch. For these reasons, as usually suggested in many-query contexts, smart reduction strategies are strongly suggested to shrink the problem dimension and to speed up the OCPs resolutions. \medskip

The Reduced Basis (RB) method \cite{Hesthaven2015, Quarteroni2015} relying, e.g., on Proper Orthogonal Decomposition (POD), may be exploited to speed up computations when solving parametrized OCPs \cite{negri2013reduced,negri2015reduced,Quarteroni2015}. Thanks to the RB method, the dimensions of state, control, and adjoint variables can be reduced through a projection onto linear subspaces of dimensions $N_y$,  $N_u$, and $N_p$, that is, we express
\begin{equation}   
\begin{alignedat}{2}
&\yh(\mus) \approx  \mathbb{V}_y \yn(\mus), \qquad \uh(\mus) \approx \mathbb{V}_u \un(\mus), \qquad
{\bf p}_h(\mus) \approx \mathbb{V}_p {\bf p}_N(\mus)
\end{alignedat}
\label{eq:POD}
\end{equation}
where the projection matrices $\mathbb{V}_y \in \R^{N_h^y \times N_y}$, $\mathbb{V}_u \in \R^{N_h^u \times N_u}$, and $\mathbb{V}_p \in \R^{N_h^p \times N_p}$ collect the basis elements spanning the low-dimensional RB subspaces. These latter are obtained through the Singular Values Decomposition (SVD) of the matrices collecting optimal state, control, and adjoint snapshots computed for suitably sampled scenarios $\{ \mus^{(1)},\ldots,\mus^{(N_s)} \} \in \mathcal{P}^{N_s}$ and, in case of time-dependent problems, on a time grid discretizing the interval $[0, T]$. In particular, RB subspaces for parametrized OCPs are the left singular vectors associated to the $N_y$, $N_u$, $N_p$ largest singular values, respectively. On the other hand, the reduced coordinates are usually obtained through a (Petrov-)Galerkin projection of each equation of \eqref{eq:nonlin_numer_optim} onto the corresponding subspace $\mathbb{V}_y$, $\mathbb{V}_u$, and $\mathbb{V}_p$. See, e.g., \citesssinline{ravindran2000reduced}{ito2001reduced}{bergmann2008optimal}{Manzoni2019}{Ravindran}{Ito and Ravindran}{Bergmann and Cordier}{Manzoni and Pagani} for some examples on the use of POD-based ROM techniques in the case of optimal control problems.   \\

Compared to an indirect approach, in which each equation appearing in the system \eqref{eq:nonlin_numer_optim} is usually reduced independently, a suitable compatibility condition must be ensured between state and adjoint reduced subspaces when a direct approach is used. In particular, if a Galerkin projection is employed to generate the ROM, an aggregated trial space $\mathbb{V}_{yp}$ made of both state and adjoint basis functions is introduced to preserve the stability of the resulting ROM. This latter would read as follows:
\begin{equation}\label{eq:KKT_compact_ROM}
\underbrace{\mathbb{V}^{\top} {\bf F}_h(\mathbb{V} {\bf q}_N (\mus))}_{\mathbf{F}_N(\mathbb{V} {\bf q}_N (\mus))} = {\bf 0}
\qquad
\mbox{with} \qquad
{\bf q}_N (\mus) = \begin{bmatrix}
{\bf y}_N(\mus) \\
{\bf u}_N(\mus) \\
{\bf p}_N(\mus) \\
\end{bmatrix},     \qquad
\mathbb{V} = \begin{bmatrix}
\mathbb{V}_{yp} & 0 & 0 \\
0 & \mathbb{V}_u & 0 \\
0 & 0 & \mathbb{V}_{yp}
\end{bmatrix}, \qquad
\mathbb{V}_{yp} = [\mathbb{V}_{y} \ \ \ \mathbb{V}_{p}]
\end{equation}
 
\medskip

\begin{remark}
In the case of a linear-quadratic OCP, the RB method would yield the following reduced order model: given $\mus$, solve 
\[
\underbrace{\mathbb{V}^{\top}  K_h(\mus) \mathbb{V}}_{K_N(\mus)} 
{\bf q}_N(\mus)  =
\underbrace{\mathbb{V}^{\top}{\bf s}_h(\mus)}_{{\bf s}_N(\mus)}.
\]
See, e.g., \citeinline{negri2013reduced, negri2015reduced}{Negri et al} for a series of applications of direct OCP solvers exploiting the RB method in the case of either advection-diffusion problems or Stokes flows, and a thorough discussion on the construction of RB spaces in those cases.
\end{remark}

\medskip

\begin{figure}
        \centering
        \includegraphics[scale = 0.45]{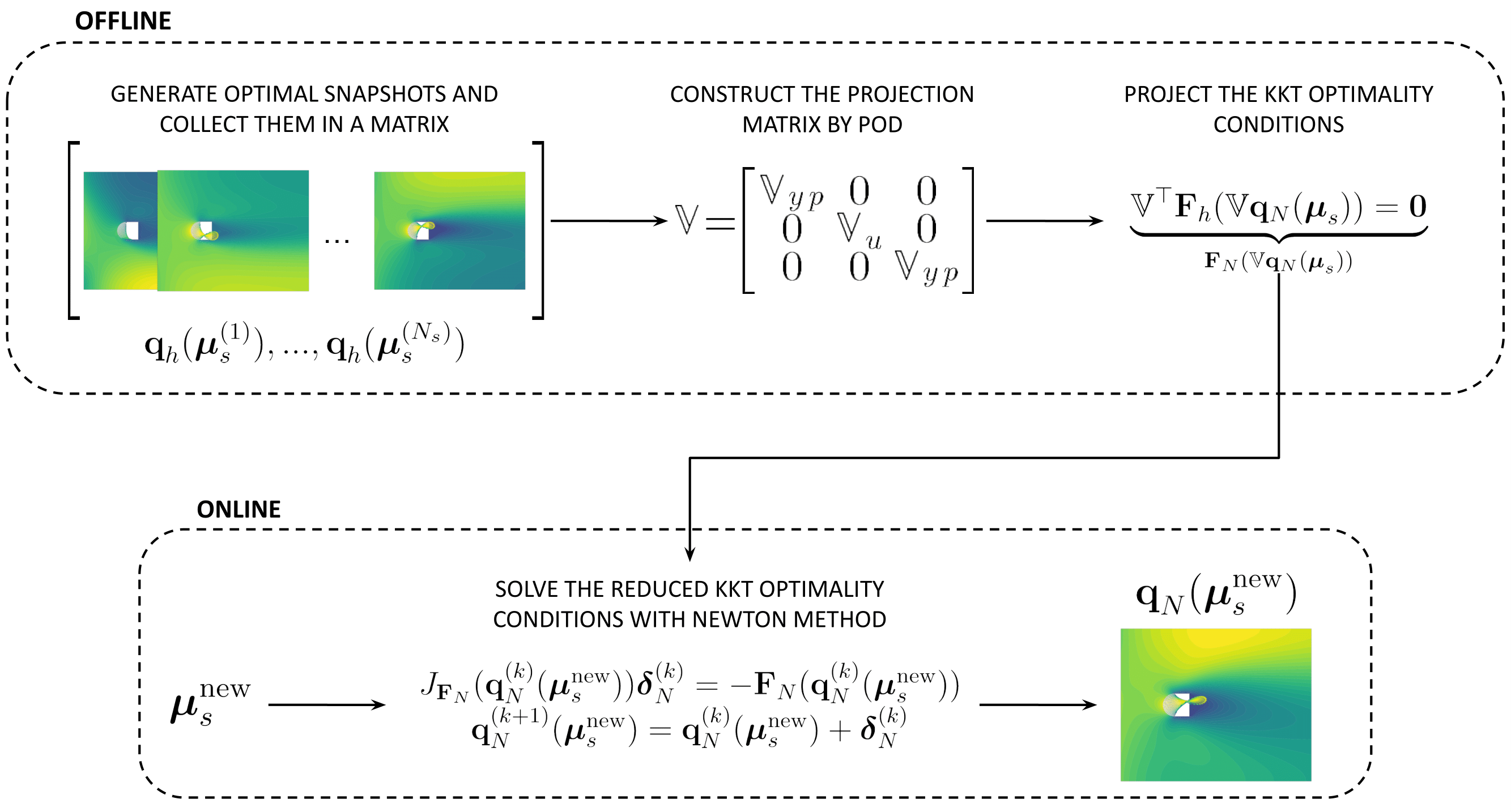}
        \caption{\textit{Test 1.1}.  Offline-online decomposition required to solve PDE-constrained optimization problems faster exploiting a reduced basis solver, in the case of steady flow control. Offline: generation of optimal state, adjoint and control snapshots for random scenarios and reduction of the first-order optimality conditions in the KKT system. Online: resolution of the reduced KKT system to compute the optimal state, adjoint and control for a new scenario of interest. In general, the reduced KKT system is solved through the Newton method with reduced Jacobian matrix $J_{\mathbf{F}_N}$ and right-hand side $\mathbf{F}_N$ properly approximated.}
        \label{fig:ParamOCP}
    \end{figure}

A direct OCP solver exploiting the RB method is summarized in Figure~\ref{fig:ParamOCP}.  For any new scenario, the optimal state, adjoint and control can be computed all-at-once by solving the reduced Karush-Kuhn-Tucker (KKT) system, that is obtained through a Galerkin projection. Thanks to the RB method, it is therefore possible to retrieve, online, the optimal solution corresponding to a new unseen scenario $\mus^{\mathrm{new}}$ faster through the resolution of a smaller KKT system.  The use of a direct approach, aiming at solving the parametrized KKT system through the RB method, can be found instead in  \citesinline{negri2013reduced,negri2015reduced}{karcher2014certified}{Negri et al}{Kärcher and Grepl}. \\

Despite several results featuring a wide range of applications, the RB method lacks of efficiency whenever aiming at the rapid control of systems in multiple scenarios. Indeed, the RB method relies on the linear superimposition of modes: hence, the ROM dimension can become moderately large to ensure accuracy in case of more complex problems, therefore avoiding a very rapid solution of the ROM problem for any new parameter instance. This might happen, e.g., if the parameter-to-solution dependence is rather involved, or in the case of complex physical behaviors of the state system (e.g., in presence of strong nonlinearities or dominated advection regimes). Moreover, the intrinsic coupling between state, adjoint, and control variables adds a further, strong difficulty, requiring to reduce all the three fields involved in the system of optimality conditions. Furthermore, a projection-based ROM is an intrusive technique since the FOM structures, such as the matrix and the right-hand side of the KKT system appearing in \eqref{eq:KKT_1}, must be directly modified and projected to obtain the ROM structures appearing in \eqref{eq:KKT_compact_ROM}. This operation is not always feasible since the arrays appearing in \eqref{eq:KKT_1} may be not directly accessible; on the other hand, the use of hyper-reduction strategies in presence of either non-affine parameter dependencies or nonlinearites is almost out of reach for OCPs given the coupled nature of these problems. Last, but not least, the inclusion of either control or state constraints is not straightforward in a ROM obtained through the RB method.  \\ 

All these difficulties led us to draw an alternative ROM strategy to tackle the rapid solution of parameterized OCPs in multiple scenarios leveraging on deep learning algorithms, focusing in particular on deep learning-based reduced order models, as detailed in the following section.

\section{Deep learning-based reduced order modeling for parametrized optimal control problems}\label{sec:OCP-DL-ROM}

This section presents a nonlinear and non-intrusive ROM strategy that allows to retrieve optimal states and controls in real-time for a given problem in multiple scenarios. Nonlinear and non-intrusive ROMs have been recently introduced -- see, e.g., \citesinline{Hesthaven2018}{Fresca2021, Fresca2022}{Hesthaven and Ubbiali}{Fresca et al} -- in order to overcome the limitations of the RB method highlighted in Section~\ref{subsec:RB}. In particular, instead of considering linear approximations relying on POD subspaces, snapshots can be compressed onto a low-dimensional space through nonlinear maps such as, e.g., autoencoders (AEs), that must be learnt. This allows to achieve the task of learning a reduced trail manifold in a much more effective way. Indeed, an autoencoder is a neural network architecture consisting of an encoder that compresses the available snapshots (to reach a very low-dimensional representation in the so-called latent space) and a decoder that recovers the full input information. The composition of the encoder and the decoder somehow approximates the identity map -- imposing that the reconstructed data are as close as possible to the original ones. Differently from POD where the encoding and decoding actions are performed through multiplications by the matrices $\mathbb{V}^{\top}$ and $\mathbb{V}$, the nonlinear activation functions within the hidden layers of the autoencoder network result in a nonlinear dimensionality reduction. For instance, \citesinline{Fresca2021}{Franco2023}{Fresca et al}{Franco et al} propose a Deep Learning-based Reduced Order Model (DL-ROM) which exploits convolutional autoencoders to take advantage of the spatial structure of PDE solutions prioritizing the dependencies of nearby values. \citeinline{Fresca2022}{Fresca and Manzoni}, instead, presents the so-called POD-DL-ROM framework where {\em (i)} the snapshots are initially compressed through POD and {\em (ii)} the resulting POD coefficients are further reduced by an autoencoder -- see also \cite{brivio2024error} for further details. The reduction combining POD and autoencoders will henceforth be referred to as POD+AE. Observe that, whenever POD is enough to retrieve an acceptable low-dimensional latent representation of the data, nonlinear reduction schemes may be avoided, as proposed by \citeinline{Hesthaven2018}{Hesthaven and Ubbiali} with the POD-NN framework.

Another key ingredient to achieve flexible and efficient ROMs is non-intrusiveness. Specifically, instead of projecting, assembling and solving a reduced system of equations as in the case of the RB method, the solutions related to new scenarios are rather retrieved by means of surrogate models. For instance, deep feed-forward neural networks may be taken into account to generate the reduced solution in the latent space starting from the input parameters $\mus$ (and, in case of time-dependent problems, the time variable $t$). For example, regarding the POD-NN framework, a neural network is exploited to approximate non-intrusively the map from $\mus$ (and possibly time) to the projected snapshots data in the POD subspace, thus avoiding the (usually, intrusive) assembling of a reduced problem. Note that the reduction step is crucial to deal with low-dimensional neural networks, allowing for lighter training and faster evaluations. By doing so, after training the solution associated with a new time instant $t^{\text{new}}$ and a new scenario $\mus^{\text{new}}$ may be computed inexpensively, through a forward pass of the parameter-to-solution map and the decoder.  

To rapidly solve a wide range of OCPs for multiple scenarios of interest, it is therefore possible to adapt and exploit the workflow of DL-ROMs introduced in \citesinline{Fresca2021}{Franco2023}{Fresca et al}{Franco et al}. Specifically, we propose a deep learning-based ROM technique in the context of parametrized OCPs, as outlined in the following offline-online decomposition.

\subsection{Offline (training) phase}

During the offline phase, the neural network architectures -- i.e., the autoencoder for the sake of dimensionality reduction, and the feedforward neural network representing the parameter-to-solution map -- must be trained. Prior to this task, a set of snapshots representing optimal (state and control) solutions for selected values of the input parameters $\mus$ must be computed exploiting a high-fidelity solver. More in detail, we perform the following steps.

\begin{itemize}
        \item[$\bullet$] We compute $N_s$ full-order optimal state and control snapshots by solving the system of first-order optimality conditions \eqref{eq:nonlin_numer_optim} for different random scenario parameters sampled in the parameter space $\mathcal{P}$, i.e.
        \begin{align*}
        & \{ \yh(\mus^{(1)}), ..., \yh(\mus^{(N_s)})\}
        \\
        & \{ \uh(\mus^{(1)}), ..., \uh(\mus^{(N_s)})\}
        \end{align*}
        where the number of snapshots $N_s$ strongly depends on several factors such as the scenario parameters dimension as well as the complexity of the parameter-to-solution map and the physics at hand. In case of time-dependent problems, for each sampled scenario, the optimal snapshots are computed on a time grid $\{t_j \}_{j=1}^{N_t}$ discretizing the interval $[0,T]$, that is
        \begin{align*}
        & \{ \yh(t_1,\mus^{(1)}), ...,\yh(t_{N_t},\mus^{(1)}),..., \yh(t_1,\mus^{(N_s)}),...,\yh(t_{N_t},\mus^{(N_s)})\}
        \\
        & \{ \uh(t_1,\mus^{(1)}), ...,\uh(t_{N_t},\mus^{(1)}),..., \uh(t_1,\mus^{(N_s)}),...,\uh(t_{N_t},\mus^{(N_s)})\}
        \end{align*}
        To generate controlled pairs, in this work we exploit high-fidelity OCP solvers based on a finite element discretization of the system of optimality conditions. Note that either direct or indirect approaches can be used to generate snapshots, in case a direct approach yields an algebraic problem hard to solve. Note also that optimal (state and control) snapshots are useful to restrict optimization procedures offline and thus to achieve real-time control capabilities in the subsequent online stage. \smallskip
        
        \item[$\bullet$] We reduce the dimensionality of optimal state snapshots and, in case of high-dimensional controls, of optimal control snapshots, too. To reach this goal either POD, an autoencoder, or a combination of the two strategies (POD+AE), can be employed. If POD is taken into account, the reduction results in the matrix multiplications between snapshots and projection matrices, that is
        \begin{equation}   
        \begin{alignedat}{2}
        &\yn =  \mathbb{V}_y^{\top} \yh, \quad 
        \un = \mathbb{V}_u^{\top} \uh
\end{alignedat}
\nonumber
\end{equation}

        \noindent When considering DL-ROMs  involving autoencoders on full-order snapshots, namely \begin{alignat*}{2}
        &\yn = \varphi_E^y(\yh); \quad &&\yh = \varphi_D^y(\yn)
        \\
        &\un = \varphi_E^u(\uh); \quad &&\uh = \varphi_D^u(\un)
        \end{alignat*}
        where $\varphi_E^y, \varphi_D^y, \varphi_E^u$ and $\varphi_D^u$ denotes the encoding and decoding networks applied respectively to state and control, the following mean squared errors are taken into account as training losses
        \begin{align*}
        & J_{AE}^y = \left\lVert \yh - \yhrec \right\rVert^2 = \left\lVert \yh - \varphi_D^y(\varphi_E^y(\yh)) \right\rVert^2
        \\
        & J_{AE}^u = \left\lVert \uh - \uhrec \right\rVert^2 = \left\lVert \uh - \varphi_D^u(\varphi_E^u(\uh)) \right\rVert^2
        \end{align*}
        where $\left\lVert\cdot\right\rVert$ stands for the Euclidean norm. Instead, in the context of POD-DL-ROM where a POD+AE reduction is performed, the autoencoder is applied on the projection of the snapshot  data onto a (potentially larger) POD space, namely 
        \begin{alignat*}{2}
        &\yn = \varphi_E^y(\mathbb{V}_y^{\top}\yh); \quad &&\yh = \mathbb{V}_y \varphi_D^y(\yn)
        \\
        &\un = \varphi_E^u(\mathbb{V}_u^{\top} \uh); \quad &&\uh = \mathbb{V}_u \varphi_D^u(\un)
        \end{alignat*}
        and the training is performed with respect to the loss functions
        \begin{align*}
        & J_{POD+AE}^y = \left\lVert \mathbb{V}_y^{\top} \yh -\mathbb{V}_y^{\top} \yhrec \right\rVert^2 = \left\lVert \mathbb{V}_y^{\top} \yh - \varphi_D^y(\varphi_E^y(\mathbb{V}_y^{\top} \yh)) \right\rVert^2
        \\
        & J_{POD+AE}^u = \left\lVert \mathbb{V}_u^{\top} \uh - \mathbb{V}_u^{\top} \uhrec \right\rVert^2 = \left\lVert \mathbb{V}_u^{\top} \uh - \varphi_D^u(\varphi_E^u(\mathbb{V}_u^{\top} \uh)) \right\rVert^2. 
        \end{align*}
        
        \item[$\bullet$] We build and train a deep feed-forward neural network $\varphi$ mapping the scenario parameters $\mus$ -- and, in case of unsteady problems, the time $t$ -- onto the latent representations of optimal state and control, that is
        \[
        \begin{bmatrix} \tildeyn \\ \tildeun \end{bmatrix} = \varphi(t, \mus).
        \]
        In particular, the loss function employed to train $\varphi$ is
        \[
        J_{\varphi} = \frac{N_y}{N_y + N_u}\left\lVert \yn - \tildeyn \right \rVert^2  + \frac{N_u}{N_y + N_u} \left\lVert \un - \tildeun \right \rVert^2 
        \]
        where the two weighting coefficients equally leverage state and control data reconstruction at the latent level. \smallskip
\end{itemize}

\begin{figure}
    \centering
    \includegraphics[scale = 0.52]{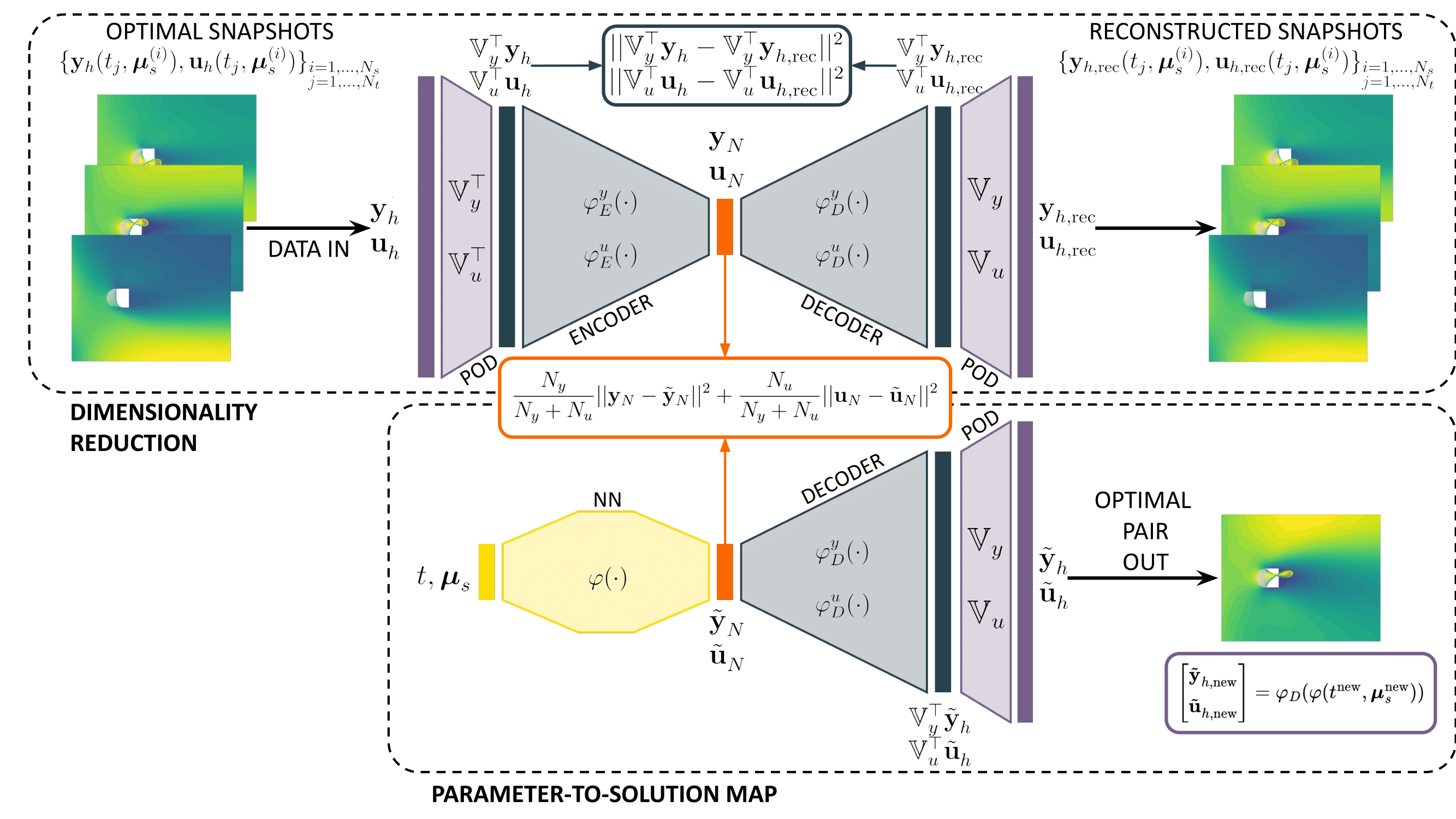}
    \caption{\textit{Test 1.2}. Deep learning-based reduced order model architecture to solve parametrized optimal control problems in real-time. Optimal state and control snapshots are generated and reduced through a combination of POD and autoencoders. Moreover, the map $\varphi$ from time and scenario parameters to reduced optimal pair is modeled through a feed-forward neural network. After properly training the networks in the offline phase, the full-order optimal pair corresponding to a new time instant $t^{\text{new}}$ and a new scenario $\mus^{\text{new}}$ is retrieved online by a forward pass through $\varphi$ and the decoders.}
    \label{fig:OCP-DL-ROM}
\end{figure}

Figure~\ref{fig:OCP-DL-ROM} visually represents the proposed architecture when taking into account a POD-DL-ROM. Note that, in general, we select the most efficient non-intrusive reduced order model -- among POD-NN, DL-ROM, and POD-DL-ROM -- for every specific OCP balancing lightweight and accuracy. In addition, according to the dynamics under investigation and the control action chosen, optimal state and control snapshots may be compressed exploiting two different reduction strategies. Moreover, $\yh$ and $\uh$ may be reduced independently, as done in the applications detailed in Section~\ref{sec:test}, or together with, e.g, a single autoencoder. The dimensionality reduction step is crucial to deal with high-dimensional state variables and distributed control fields, as of interest in the applications presented in Section~\ref{sec:test}: indeed, it allows to consider parameter-to-solution maps with compressed output sizes that are lighter to train and faster to evaluate online.

\begin{remark}
Note that, despite a set of solutions of the adjoint problem are generated while solving the system of optimality conditions during the snapshots' calculation, these data are discarded and no reduction is operated on the adjoint field, in the proposed workflow. This also represents a major difference compared to projection-based ROMs for optimal control problems, that must include the evaluation of the adjoint solution, thus yielding higher complexity and costs. Nevertheless, in case the adjoint solution is also required as a specific requirement of the application at hand, the workflow discussed so far can be easily adapted to generate the adjoint field, too.
\end{remark}

\subsection{Online phase}

During the online phase, once the proposed architecture has been trained, it is possible to evaluate the optimal (state and control) solution for new scenarios. In particular, we can  retrieve the optimal pair associated with a new scenario $\mus^{\mathrm{new}}$ and, in case of time-dependent problems, a new time instant $t^{\mathrm{new}}$ through a forward pass of the trained neural networks, that is
        \[
        \begin{bmatrix} \tilde{\mathbf{y}}_{h,\text{new}} \\ \tilde{\mathbf{u}}_{h,\text{new}}\end{bmatrix} = \varphi_D\left(\begin{bmatrix} \tilde{\mathbf{y}}_{N,\text{new}} \\ \tilde{\mathbf{u}}_{N,\text{new}} \end{bmatrix} \right) = \varphi_D(\varphi(t^{\mathrm{new}}, \mus^{\mathrm{new}}))
        \]
        where, for the sake of compactness, $\varphi_D$ combines the state and control decoders. Specifically, depending on the reduction technique taken into account, $\varphi_D$ results in 
        \[
        \varphi_D\left(\begin{bmatrix} \tilde{\mathbf{y}}_{N,\text{new}} \\ \tilde{\mathbf{u}}_{N,\text{new}} \end{bmatrix} \right) = \begin{cases} 
        \begin{bmatrix} \mathbb{V}_y & 0 \\ 0 & \mathbb{V}_u \end{bmatrix} \begin{bmatrix} \tilde{\mathbf{y}}_{N,\text{new}} \\ \tilde{\mathbf{u}}_{N,\text{new}} \end{bmatrix} \quad & \text{in case a POD-NN is used} \vspace{0.2 cm}
        \\
        \begin{bmatrix} \varphi_D^y \left( \tilde{\mathbf{y}}_{N,\text{new}} \right) \\ \varphi_D^u \left( \tilde{\mathbf{u}}_{N,\text{new}} \right) \end{bmatrix} & \text{in case a DL-ROM is used} \vspace{0.2 cm}
        \\
        \begin{bmatrix} \mathbb{V}_y & 0 \\ 0 & \mathbb{V}_u \end{bmatrix} \begin{bmatrix} \varphi_D^y \left( \tilde{\mathbf{y}}_{N,\text{new}} \right) \\ \varphi_D^u \left( \tilde{\mathbf{u}}_{N,\text{new}} \right) \end{bmatrix} & \text{in case a POD-DL-ROM is used}. 
        \end{cases} \medskip
        \]
The proposed architecture is extremely flexible thanks to its non-intrusive data-driven nature and the independence of the surrogate models on the particular problem faced. Indeed, it can be applied to a wide range of parametrized OCPs, as demonstrated throughout Section~\ref{sec:test}, and can be extended to deal with different data sources, such as sensors or recordings.

\section{Numerical results}\label{sec:test}

This section is devoted to the numerical results obtained with the proposed architecture when tackling three distinct parametric OCPs in two spatial dimensions. Section~\ref{subsec:NS} focuses on nonlinear OCPs involving the minimization of the energy dissipated by a fluid flow passing an obstacle in a channel considering both steady and time-dependent Navier-Stokes equations to model its dynamics. Section~\ref{subsec:cooling}, instead, takes into account an active thermal cooling problem where the temperature of an object has to be kept equal to a reference value. To perform error analysis and fairly assess the performance of the proposed machinery, we consider a test dataset $\{ \yh(t_j, \mus^{(i)}), \uh(t_j, \mus^{(i)})\}_{(i,j) \in I_{test}}$ where $I_{test} \subset \{1,...,N_s \} \times \{1,...,N_t \}$ and $|I_{test}| \approx 0.2 N_t N_s$. In particular, the following mean relative errors are considered for evaluation purposes:
\[
\varepsilon^y_{\mathrm{rel,}L^2} = \dfrac{1}{|I_{test}|} \sum_{(i,j) \in I_{test}} 
 \sqrt{\dfrac{ \int_{\Omega} \lVert y_h(t_j, \mus^{(i)}) - \tilde{y}_h (t_j, \mus^{(i)})\rVert^2 d\Omega }{ \int_{\Omega} \lVert y_h(t_j, \mus^{(i)})\rVert^2 d\Omega }}\myspace
 \varepsilon^u_{\mathrm{rel,}L^2} = \dfrac{1}{|I_{test}|} \sum_{(i,j) \in I_{test}} 
 \sqrt{\dfrac{ \int_{\Omega} \lVert u_h(t_j, \mus^{(i)}) - \tilde{u}_h (t_j,\mus^{(i)})\rVert^2 d\Omega }{ \int_{\Omega} \lVert u_h(t_j,\mus^{(i)})\rVert^2 d\Omega }}
\]
where $y_h$ and $u_h$ are the ground truth finite element functions associated with the coefficients $\yh$ and $\uh$. Instead, $\tilde{y}_h$ and $\tilde{u}_h$ are the corresponding approximations associated with $\yhrec$ and $\uhrec$ when evaluating the ability of the chosen reduction technique, while they are linked to $\tildeyh$ and $\tildeuh$ when interested in the prediction errors committed by the proposed architecture. Notice that the Euclidean norm $\lVert \cdot \rVert$ is useful to take into account both scalar and vector functions.

\subsection{Flow control problems}
\label{subsec:NS}

\begin{figure}
    \centering
    \includegraphics[height = 0.35\linewidth]{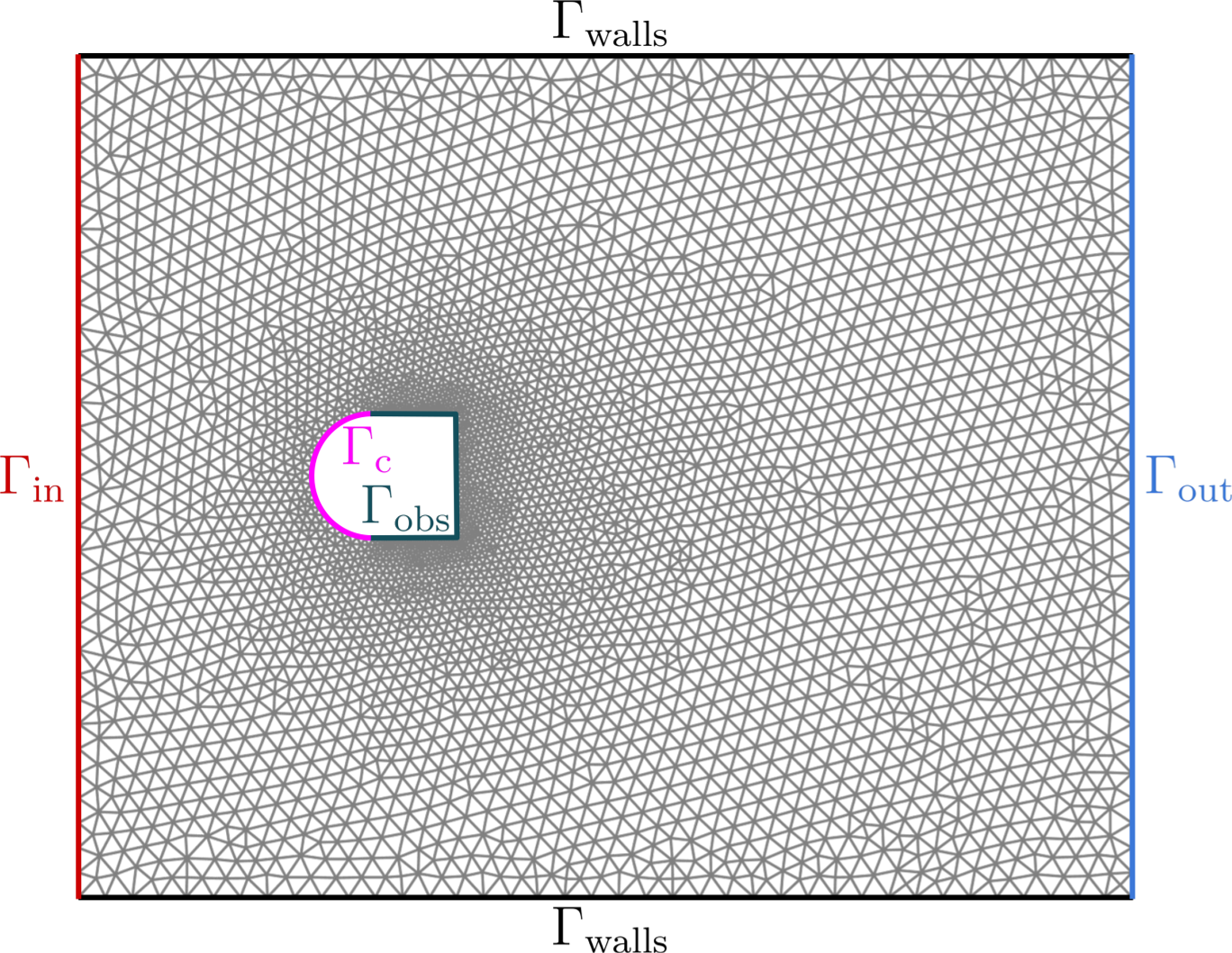}
   \caption{\textit{Test 1.} Flow control. Mesh exploited to generate high-fidelity snapshots both in steady and time-dependent flow control test cases, along with the boundaries considered in Navier-Stokes equations. In particular, $\Gamma_{\mathrm{in}}$ in red is the inflow boundary where the fluid enters the channel, $\Gamma_{\mathrm{out}}$ in blue is the outflow, $\Gamma_{\mathrm{walls}}$ in black are the walls bounding the channel from the top and the bottom, $\Gamma_{\mathrm{obs}}$ in green and the control region $\Gamma_c$ in magenta are the two portions of the obstacle boundary.}
 \label{fig:NS_mesh}
\end{figure}

In this section, we test the proposed real-time OCP solver in the case of nonlinear Navier-Stokes equations, which are resumed in the following system:

\begin{equation}
    \begin{cases}
       \dfrac{\partial \mathbf{v}(\mathbf{x},t)}{\partial t} -\nu \Delta \mathbf{v}(\mathbf{x},t) + (\mathbf{v}(\mathbf{x},t) \cdot \nabla) \mathbf{v}(\mathbf{x},t) + \nabla p(\mathbf{x},t) = 0  \qquad &\mathrm{in} \ \Omega \times (0,T) \\
       \mathrm{div }\  \mathbf{v}(\mathbf{x},t) = 0  \qquad &\mathrm{in} \ \Omega \times (0,T) \\
       \mathbf{v}(\mathbf{x},t) = \mathbf{u}(\mathbf{x},t)  \qquad &\mathrm{on} \ \Gamma_{\mathrm{c}} \times (0,T]\\
       \mathbf{v}(\mathbf{x},t) = \mathbf{0}  \qquad &\mathrm{on} \ \Gamma_{\mathrm{obs}} \times (0,T]\\
       \mathbf{v}(\mathbf{x},t) = \mathbf{v}_{\text{in}}(\mus)  \qquad &\mathrm{on} \ \Gamma_{\mathrm{in}} \times (0,T]\\
       \mathbf{v}(\mathbf{x},t) \cdot \mathbf{n}(\mathbf{x}) = 0  \qquad &\mathrm{on} \ \Gamma_{\mathrm{walls}} \times (0,T] \\
       (\nu \nabla \mathbf{v}(\mathbf{x},t) - p(\mathbf{x},t))\mathbf{n}(\mathbf{x}) \cdot \mathbf{t}(\mathbf{x}) = 0  \qquad &\mathrm{on} \ \Gamma_{\mathrm{walls}} \times (0,T] \\
       (\nu \nabla \mathbf{v}(\mathbf{x},t) - p(\mathbf{x},t))\mathbf{n}(\mathbf{x}) = 0  \qquad &\mathrm{on} \ \Gamma_{\mathrm{out}} \times (0,T]\\
       \mathbf{v}(\mathbf{x},0) = \mathbf{0} &\mathrm{in} \ \Omega \times \{t = 0\}.
\end{cases}
\label{eq:NS}
\end{equation}
\noindent In the equations above, $\nu$ denotes the kinematic viscosity, while $\mathbf{n}$ and $\mathbf{t}$ are the unit vectors normal and tangential to the domain boundary, respectively. We consider a rectangular channel $\Omega = [0,25] \times [0,20]$ with an obstacle inside, inspired to configurations already explored in the literature several times \cite{Dedè2007,negri2015reduced,Manzoni2021}. The domain discretization in Figure~\ref{fig:NS_mesh} is generated through \texttt{gmsh} utilities \citep{gmsh} and consists of $3716$ vertices and $7203$ elements. The state variable in Equation~\eqref{eq:NS} includes both the velocity vector $\mathbf{v}(\mathbf{x},t)$ ($m s^{-1}$) and the pressure $p(\mathbf{x},t)$ ($Pa$). However, for simplicity, only the velocity is considered in the following, disregarding the pressure field -- although also this latter field could be retrieved by means of a slight adaptation of the proposed framework. 

When dealing with the full-order model, the state is discretized with $\mathbb{P}_2-\mathbb{P}_1$ finite elements, resulting in $N_h^{\mathbf{v}} = 29270$ degrees of freedom for the velocity. The control action is instead represented by the time-dependent Dirichlet velocity datum $\mathbf{u}$ on the rounded edge of the obstacle boundary $\Gamma_c$, which can be interpreted as fluid injection or absorption. $\mathbf{u}$ is thus a boundary control depending on $N_h^{\mathbf{u}} = 98$ degrees of freedom at the high-fidelity level when discretized through FEM. As far as boundary and initial conditions are concerned, the following are taken into account: homogeneous Dirichlet conditions on the right portion of the boundary $\Gamma_{\mathrm{obs}}$, homogeneous Neumann boundary conditions on the outflow $\Gamma_{\mathrm{out}}$, free-slip boundary conditions on the upper and lower walls $\Gamma_{\mathrm{walls}}$, and homogeneous initial conditions for $t=0$ seconds. On the inflow boundary $\Gamma_{\mathrm{in}}$, Dirichlet boundary conditions are imposed with velocity datum given by
\[
\mathbf{v}_{\mathrm{in}} = \left\lVert \mathbf{v} \right\rVert_{\mathrm{in}}[\cos(\alpha_{\mathrm{in}}), 0.01 x_2 (20 - x_2) \sin(\alpha_{\mathrm{in}})]^{\top}
\]
where $\alpha_{\mathrm{in}}$ (radians) denotes the angle of attack and $\left\lVert \mathbf{v} \right\rVert_{\mathrm{in}}$ ($m s^{-1}$) the inflow intensity.

In this setting, we aim to minimize the energy dissipated by the fluid, that can be expressed as the $L^2$ norm of the velocity gradient. The loss function is enriched by two regularizing terms dealing with the $L^2$ norm of $\mathbf{u}$ and its gradient, in order to stabilize the control problem and to avoid extremely expensive and energetic optimal control strategies. In the end, the selected loss function is 
\[
J(\mathbf{v}(\mathbf{x},t),\mathbf{u}(\mathbf{x},t)) = \int_0^T \int_{\Omega} \left \lVert \nabla \mathbf{v}(\mathbf{x},t) \right\rVert^2 d\Omega dt + \beta \int_0^T \int_{\Gamma_c} \left \lVert\mathbf{u}(\mathbf{x},t)\right\rVert^2 d\Gamma_c dt + \beta_g \int_0^T \int_{\Gamma_c}\left \lVert \nabla \mathbf{u}(\mathbf{x},t)\right \rVert^2 d\Gamma_c dt
\]
where the time integrals are considered only for time-dependent settings and the coefficients $\beta$ and $\beta_g$ are introduced to balance the integrals in the cost functional that may have different magnitudes. In particular, $\beta = \beta_g = 10^{-2}$ in order to prioritize the minimization of energy dissipation. Note that, after discretizing Equation~\eqref{eq:NS} through FEM, it is possible to obtain the discretized cost functional $J_h$ introduced in Section~\ref{sec:OCP}.

\subsubsection{Flow control for steady Navier-Stokes equations}
\label{subsec:steadyNS}

\begin{figure}
    \centering
    \subfloat{\includegraphics[height = 0.3\linewidth]{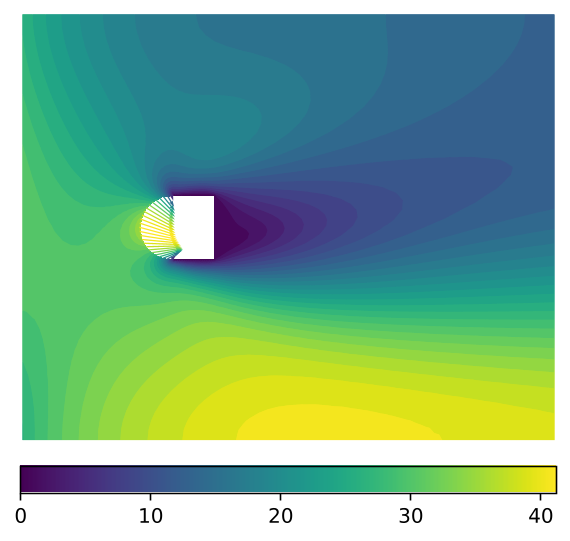}}    \quad \subfloat{\includegraphics[height = 0.3\linewidth]{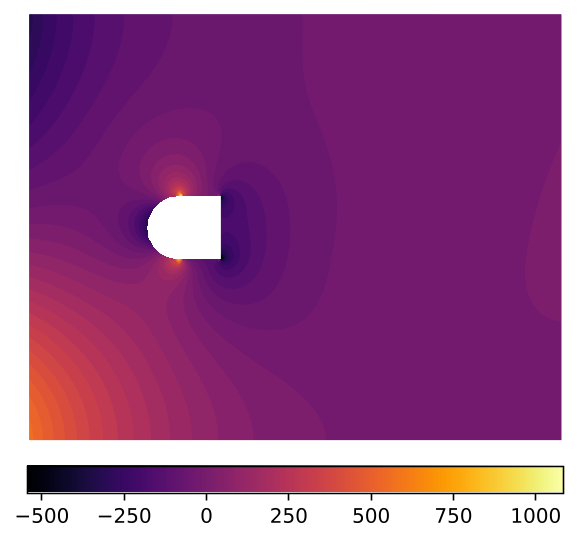}}
    
   \caption{\textit{Test 1.1}. Steady  flow control. Optimal control, velocity and pressure obtained for $\left\lVert\mathbf{v}_{\mathrm{in}}\right\rVert = 31.16 m s^{-1}$ and $\alpha_{\mathrm{in}} = -0.53$ radians through the high-fidelity OCP solver. The velocity on $\Omega$ is depicted through a scalar field with colours corresponding to its norm, while the control on $\Gamma_c$ is represented through a vector field.}
 \label{fig:steadyNS_setting}
\end{figure}

We first consider the steady Navier-Stokes equations, reported here for the sake of completeness
\begin{equation}
    \begin{cases}
        -\nu \Delta \mathbf{v}(\mathbf{x}) + (\mathbf{v}(\mathbf{x}) \cdot \nabla) \mathbf{v}(\mathbf{x}) + \nabla p(\mathbf{x}) = 0  \qquad &\mathrm{in} \ \Omega \\
       \mathrm{div }\  \mathbf{v}(\mathbf{x}) = 0  \qquad &\mathrm{in} \ \Omega \\
       \mathbf{v}(\mathbf{x}) = \mathbf{u}(\mathbf{x})  \qquad &\mathrm{on} \ \Gamma_{\mathrm{c}} \\
       \mathbf{v}(\mathbf{x}) = \mathbf{0}  \qquad &\mathrm{on} \ \Gamma_{\mathrm{obs}} \\
       \mathbf{v}(\mathbf{x}) = \mathbf{v}_{\text{in}}(\mus)  \qquad &\mathrm{on} \ \Gamma_{\mathrm{in}} \\
       \mathbf{v}(\mathbf{x}) \cdot \mathbf{n}(\mathbf{x}) = 0  \qquad &\mathrm{on} \ \Gamma_{\mathrm{walls}} \\
       (\nu \nabla \mathbf{v}(\mathbf{x}) - p(\mathbf{x}))\mathbf{n}(\mathbf{x}) \cdot \mathbf{t}(\mathbf{x}) = 0  \qquad &\mathrm{on} \ \Gamma_{\mathrm{walls}} \\
       (\nu \nabla \mathbf{v}(\mathbf{x}) - p(\mathbf{x}))\mathbf{n}(\mathbf{x}) = 0  \qquad &\mathrm{on} \ \Gamma_{\mathrm{out}}
\end{cases}
\nonumber
\end{equation}
\noindent with $\nu = 10 m^2 s^{-1}$. The two parameters characterizing the inflow are regarded as scenario parameters, that is $\mus = [\alpha_{\mathrm{in}}, \left\lVert\mathbf{v}\right\rVert_{\mathrm{in}}]^{\top}$. Following the procedure detailed in Section~\ref{sec:OCP-DL-ROM}, we generate $N_s = 50$ optimal pairs related to random inflow intensities and angles of attack uniformly sampled in the intervals $(10.0, 130.0) m s^{-1}$ and $(-1.0, 1.0)$ radians, respectively. The Reynolds number is thus very low and varies between $3$ and $39$. The forward resolution of Equation~\eqref{eq:NS} is performed through the Newton method implemented in \texttt{fenics} \citep{fenics}, while optimal state and control snapshots are calculated with an OCP solver based on finite elements and the adjoint method implemented in \texttt{dolfin-adjoint} \citep{Mitusch2019}. In particular, the high-fidelity solver takes, on average, $11$ minutes to retrieve every snapshot. Out of all simulated data, $40$ optimal pairs are exploited for dimensionality reduction and to train neural networks (training set), while the remaining ones are used for testing purposes (test set). The optimal snapshot generated for $\left\lVert\mathbf{v}_{\mathrm{in}}\right\rVert = 31.16 m s^{-1}$ and $\alpha_{\mathrm{in}} = -0.53$ radians is shown in Figure~\ref{fig:steadyNS_setting}. \\

\begin{figure}
\centering
\subfloat{\includegraphics[height = 0.3\linewidth]{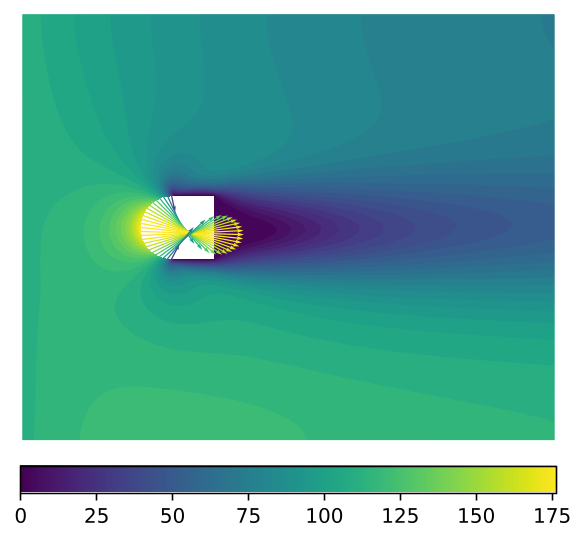}}
\subfloat{\includegraphics[height = 0.3\linewidth]{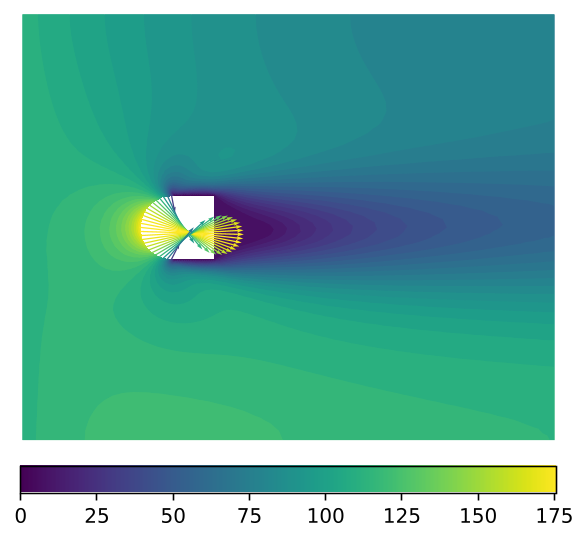}}
\subfloat{\includegraphics[height = 0.3\linewidth]{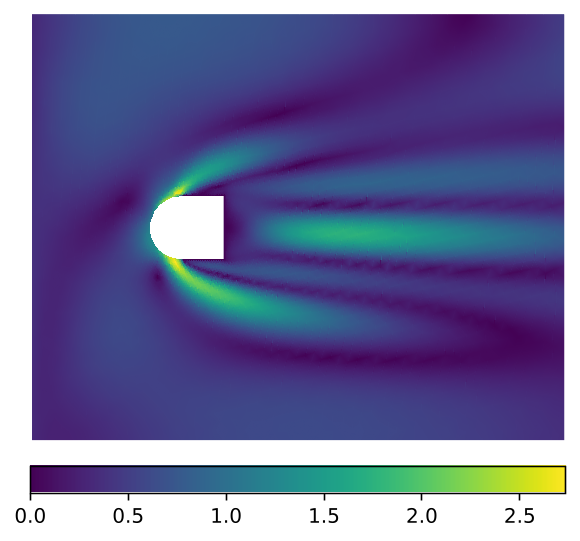}}

\hphantom{\subfloat{\includegraphics[height = 0.3\linewidth]{Images/steadyNS/Ground_truth_modin=114_alphain=-0.17.png}}}\subfloat{\includegraphics[height = 0.3\linewidth]{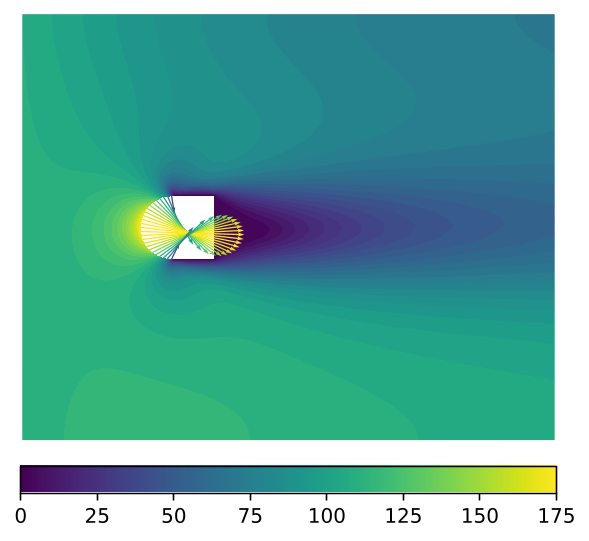}}
\subfloat{\includegraphics[height = 0.3\linewidth]{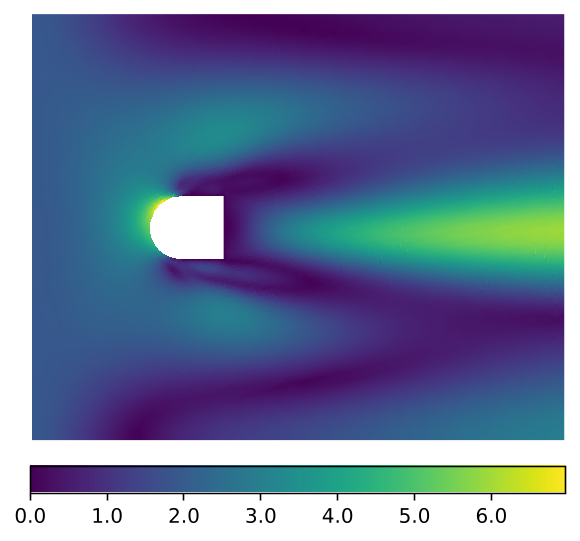}}

\caption{\textit{Test 1.1}. Steady flow control. First row: high-fidelity optimal snapshot, POD reconstruction and reconstruction error corresponding to the test scenario parameters $\left\lVert\mathbf{v}_{\mathrm{in}}\right\rVert = 114.79m s^{-1}$ and $\alpha_{\mathrm{in}} = -0.17$ radians. Second row: POD-NN prediction and prediction error corresponding to the test scenario parameters $\left\lVert\mathbf{v}_{\mathrm{in}}\right\rVert = 114.79m s^{-1}$ and $\alpha_{\mathrm{in}} = -0.17$ radians. The velocity on $\Omega$ is depicted through a scalar field with colors
corresponding to its norm, while the control on $\Gamma_c$ is represented through a vector field.}
\label{fig:steadyNS_1}
\end{figure}

\begin{figure}
\centering
\subfloat{\includegraphics[height = 0.3\linewidth]{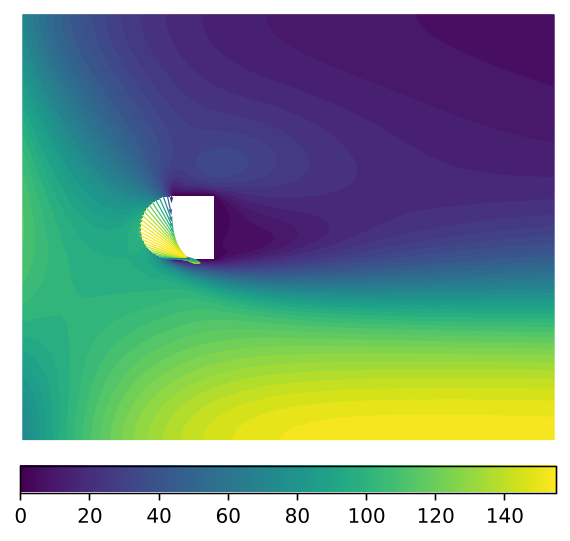}}
\subfloat{\includegraphics[height = 0.3\linewidth]{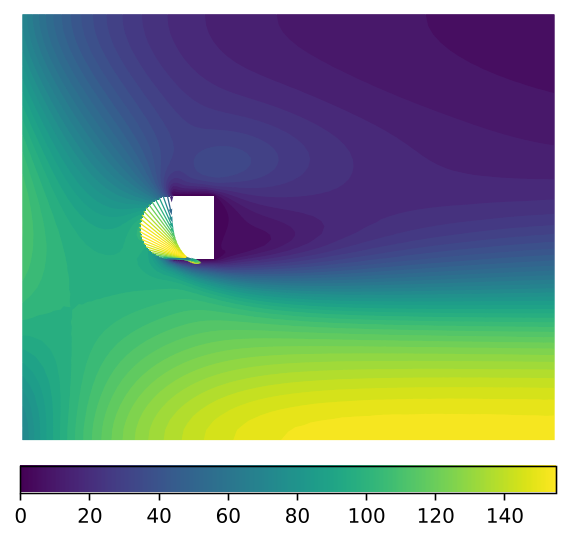}}
\subfloat{\includegraphics[height = 0.3\linewidth]{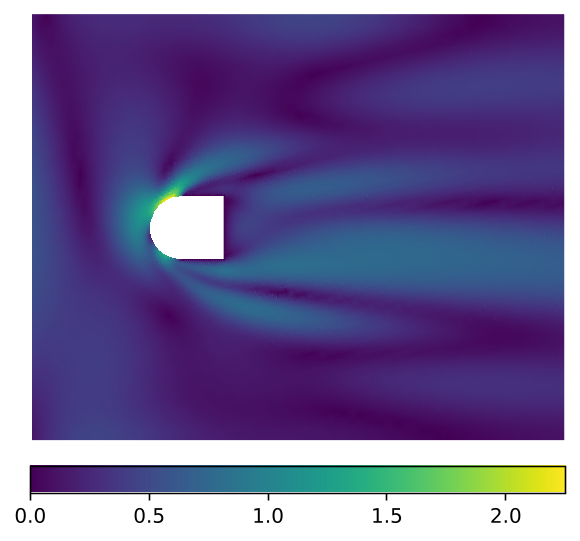}}

\hphantom{\subfloat{\includegraphics[height = 0.3\linewidth]{Images/steadyNS/Ground_truth_modin=122.7_alphain=-0.89.png}}}\subfloat{\includegraphics[height = 0.3\linewidth]{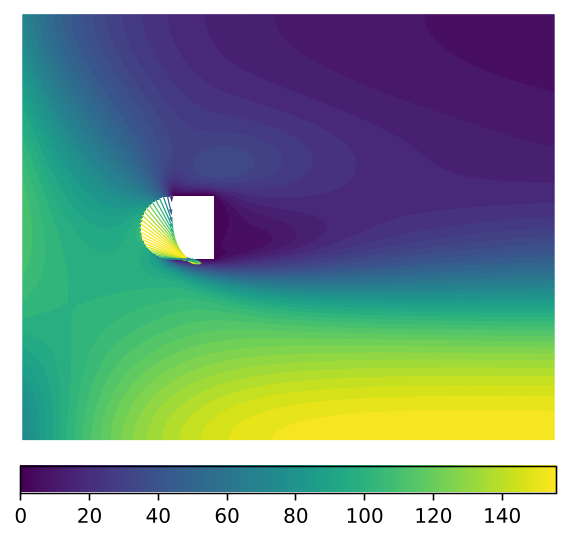}}
\subfloat{\includegraphics[height = 0.3\linewidth]{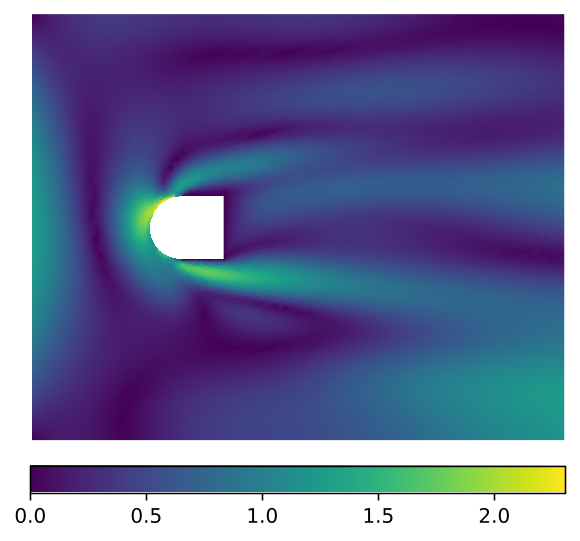}}

  \caption{\textit{Test 1.1}. Steady flow control. First row: high-fidelity optimal snapshot, POD reconstruction and reconstruction error corresponding to the test scenario parameters $\left\lVert\mathbf{v}_{\mathrm{in}}\right\rVert = 112.71m s^{-1}$ and $\alpha_{\mathrm{in}} = -0.88$ radians. Second row: POD-NN prediction and prediction error corresponding to the test scenario parameters $\left\lVert\mathbf{v}_{\mathrm{in}}\right\rVert = 112.71m s^{-1}$ and $\alpha_{\mathrm{in}} = -0.88$ radians. The velocity on $\Omega$ is depicted through a scalar field with colours corresponding to its norm, while the control on $\Gamma_c$ is represented through a vector field.}
\label{fig:steadyNS_2}
\end{figure}

To build a faster but still reliable OCP solver following the architecture presented in Section~\ref{sec:OCP-DL-ROM}, we start reducing the snapshots dimensionality through POD. In particular, when reducing vector fields, POD is applied component-wise, i.e. the $x_1$ and $x_2$ components of the velocity are reduced separately, while considering the same number of modes for the two components. In this test case, nonlinear reduction strategies based on autoencoders are not necessary thanks to a fast singular value decay that allows to achieve acceptable reduction errors with few POD modes. Specifically, considering $N_{\mathbf{u}} = 6$ and $N_{\mathbf{v}} = 10$, the $L^2$ mean relative error committed when trying to reconstruct test data is $0.73 \%$ on the state and $0.70 \%$ on the control. Top rows of Figure~\ref{fig:steadyNS_1} and Figure~\ref{fig:steadyNS_2} display the ground truth and the reconstructed velocity fields, that are the snapshots compressed and reconstructed through POD, related to two different choices of scenario parameters in the test set.

To rapidly retrieve the optimal state and control values for different inflow velocities and angles of attack, here we must learn the map from the scenario parameters onto the POD coefficients of the optimal pairs, thus following a POD-NN strategy. This is done considering a feed-forward neural network $\varphi$ having $3$ hidden layers with, respectively, $50, 100, 50$ neurons and exploiting leaky Relu as activation function. After building and training the neural network in $26.84$ seconds through \texttt{pytorch} utilities taking into account the L-BFGS optimization algorithm, the predicted test data are accurate up to a $L^2$ mean relative error equal to $1.53\%$ for the optimal state and $1.81\%$ for the optimal control. The bottom rows of Figure~\ref{fig:steadyNS_1} and Figure~\ref{fig:steadyNS_2} display the ground truth optimal solutions related to two different choices of scenario parameters in the test set and the corresponding POD-NN predictions given by
\[
\begin{bmatrix} \tilde{\mathbf{v}}_h \\ \tildeuh \end{bmatrix} = \begin{bmatrix} \mathbb{V}_y & 0 \\ 0 & \mathbb{V}_u \end{bmatrix} \varphi(\alpha_{\mathrm{in}}, \left\lVert\mathbf{v}\right\rVert_{\mathrm{in}}
)
\]        
 
To emphasize the importance of dimensionality reduction when tackling high-dimensional problems, we fit a neural network mapping the scenario parameters directly into the high-fidelity optimal state and control variables, resulting in an output dimension equal to $N_h^{\mathbf{v}} + N_h^{\mathbf{u}} = 29368$. To achieve a $L^2$ mean relative error of about $4\%$ for both state and control, which is acceptable but higher than the one obtained with POD-NN, a neural network involving $135$ times more weights,  resulting in a training time $400$ times longer, with respect to the one required for POD-NN shall be required.

\begin{figure}
    \centering
    \subfloat{\includegraphics[scale = 0.5]{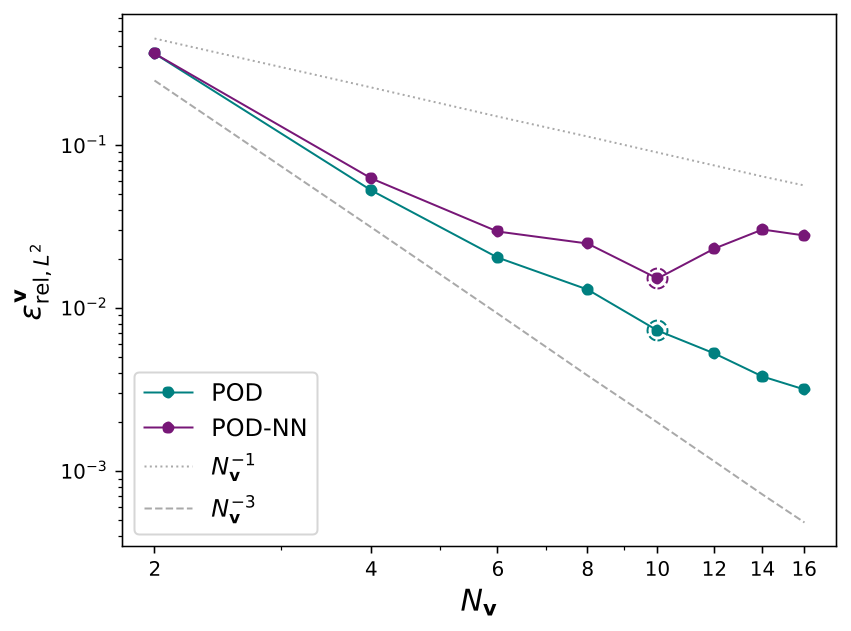}}
    \quad
    \subfloat{
    \includegraphics[scale = 0.5]{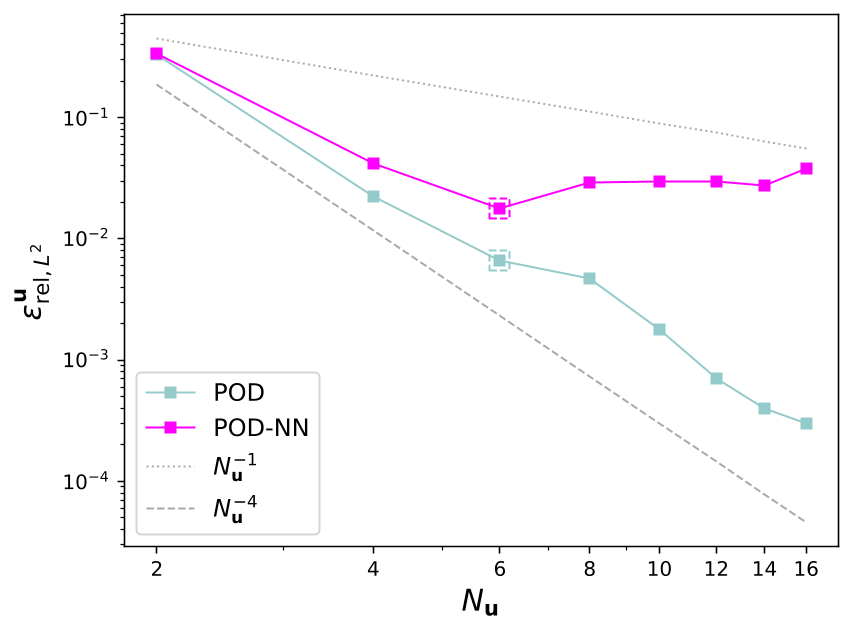}}
   
    \caption{\textit{Test 1.1}. Steady flow control. Left: $L^2$ mean relative error on $\mathbf{v}$ committed by POD (reconstruction error) and POD-NN (prediction error) for different latent dimensions $N_{\mathbf{v}}$. Right: $L^2$ mean relative error on $\mathbf{u}$ committed by POD (reconstruction error) and POD-NN (prediction error) for different latent dimensions $N_{\mathbf{u}}$. The latent dimensions selected in the numerical example are highlighted with dashed markers.}
   \label{fig:steadyNS_error_analysis1}
\end{figure}

Figure~\ref{fig:steadyNS_error_analysis1} proposes a comparison between the reconstruction errors offered by POD and the prediction errors committed by POD-NN for different latent dimensions, while keeping fixed the other hyperparameters such as the $\varphi$ architecture and the training set dimension. The decay of the POD reconstruction errors on test data, both for the state and the control, are superlinear. The POD-NN errors on test data are higher than those provided by POD due to the approximation given by the neural network $\varphi$, while keeping similar superlinear decays with respect to the latent dimension both for the control and for the state approximation. Furthermore, it is possible to assess a saturation in the POD-NN performance for large latent dimensions.

\begin{figure}
    \centering
    \subfloat{\includegraphics[height = 0.3\linewidth]{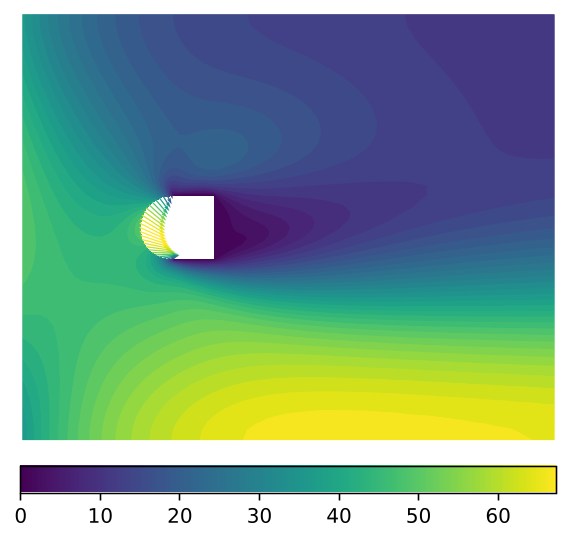}} \quad 
    \subfloat{\includegraphics[height = 0.3\linewidth]{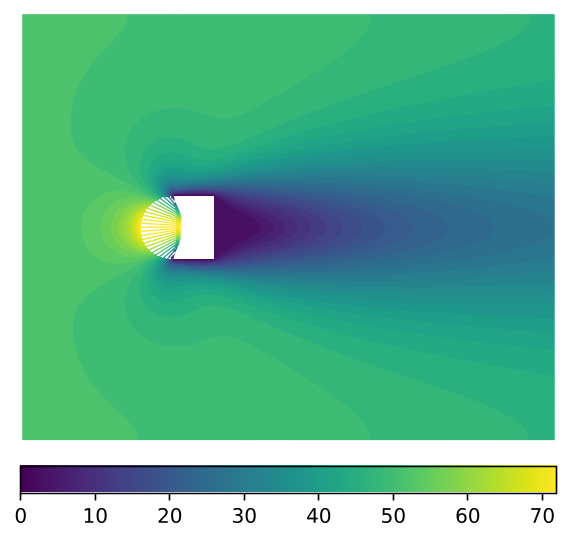}} \quad \subfloat{\includegraphics[height = 0.3\linewidth]{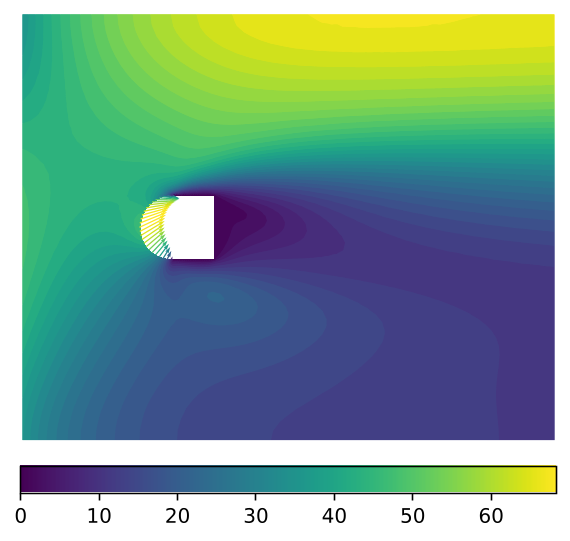}}

  \caption{\textit{Test 1.1}. Steady flow control. Optimal state and control provided by POD-NN corresponding to the scenario parameters $\left\lVert\mathbf{v}_{\mathrm{in}}\right\rVert = 100m s^{-1}$ and $\alpha_{\mathrm{in}} = -0.75, 0, 0.75$ radians. The velocity on $\Omega$ is depicted through a scalar field with colours corresponding to its norm, while the control on $\Gamma_c$ is represented through a vector field.}
\label{fig:steadyNS_results}
\end{figure}

The POD-NN architecture is then exploited to solve the parametric OCP for new values of scenario parameters unseen during the data generation and training phases. In particular, we consider $\left\lVert\mathbf{v}_{\mathrm{in}}\right\rVert = 100m s^{-1}$ and several angles of attack in a grid spanning $(-1, 1)$ radians. Figure~\ref{fig:steadyNS_results} shows the results obtained for three different angles of attack considered. Each result is computed online in less than $0.002$ seconds thanks to a forward pass through the neural network $\varphi$ and the POD decoding in order to predict full-order optimal pairs.

\begin{figure}
    \centering
    \subfloat{\includegraphics[scale = 0.5]{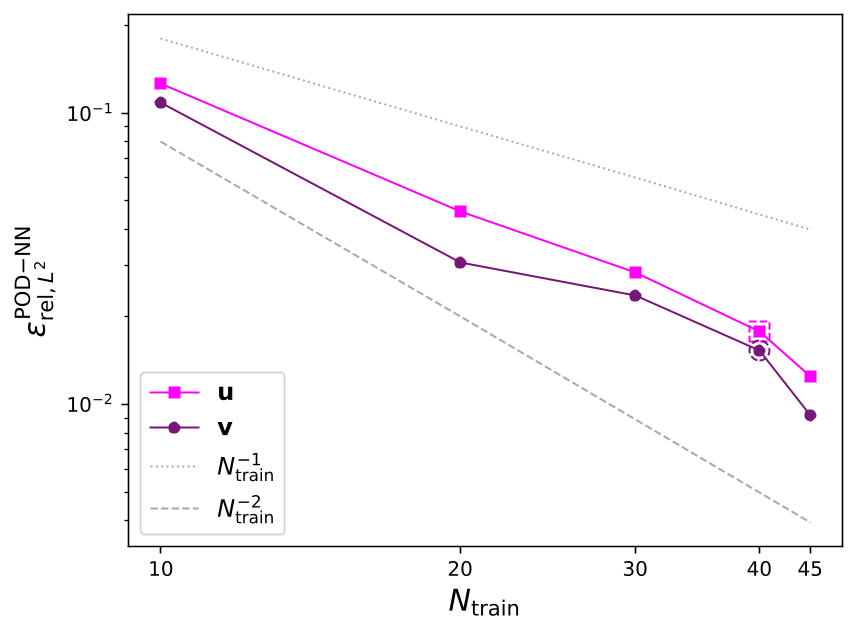}}
    \quad
    \subfloat{\includegraphics[scale = 0.5]{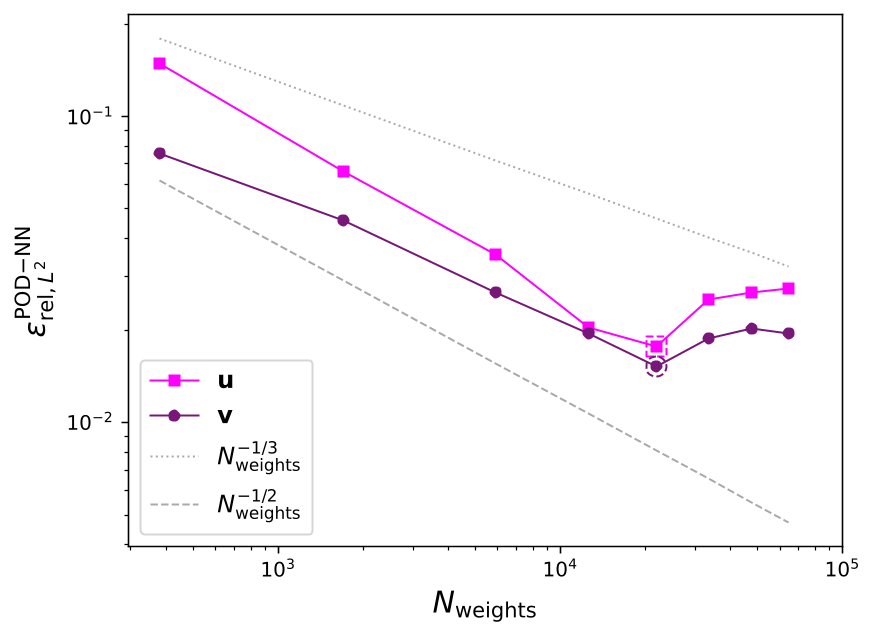}}
    
    \caption{\textit{Test 1.1}. Steady flow control. Left: $L^2$ mean relative test error decay with respect to training set dimension $N_{\mathrm{train}}$. Right: $L^2$ mean relative test error decay with respect to the number of $\varphi$ weights $N_{\mathrm{weights}}$. The hyperparameters selected in the numerical example are highlighted with dashed markers.}
   \label{fig:steadyNS_error_analysis2}
\end{figure}

Figure~\ref{fig:steadyNS_error_analysis2} presents, instead, two different analyses that are helpful in understanding the benefits and behavior of the proposed method. In particular, the left panel investigates the trend of $L^2$ mean relative test error for different training set dimensions $N_{\mathrm{train}}$ keeping fixed the other hyperparameters, while a similar analysis is repeated in the right panel for different $\varphi$ complexities, i.e. for different number of neural network weights $N_{\mathrm{weights}}$. The superlinear decay obtained in the first case suggests that an increase in the number of training data has an higher impact on the goodness of fit of the proposed model rather than an increase of the number of $\varphi$ weights, which entails a sublinear decay.

\subsubsection{Flow control for time-dependent Navier-Stokes equations}
\label{subsec:unsteadyNS}

\begin{figure}
    \centering
    \subfloat{\includegraphics[height = 0.3\linewidth]{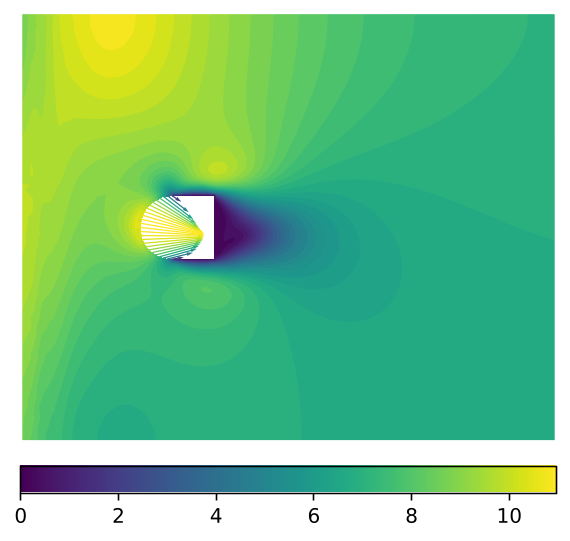}}    \quad \subfloat{\includegraphics[height = 0.3\linewidth]{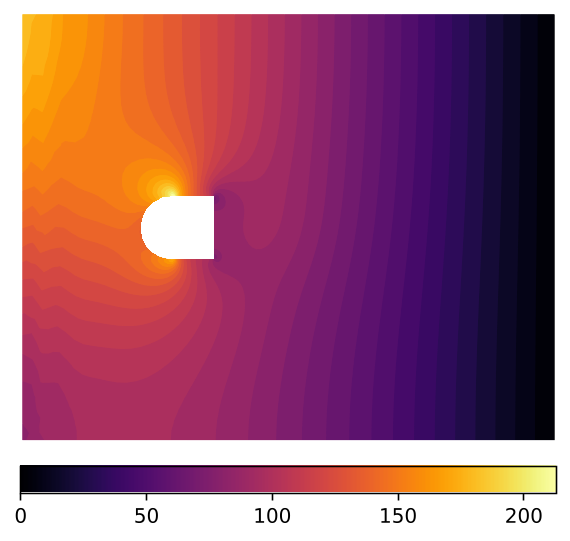}}
    
   \caption{\textit{Test 1.2}. Time-dependent flow control. Optimal control, velocity and pressure obtained for $\alpha_{\mathrm{in}} = 0.45$ radians and $t = 0.5$ seconds through a high-fidelity full-order OCP solver. The velocity on $\Omega$ is depicted through a scalar field with colours corresponding to its norm, while the control on $\Gamma_c$ is represented through a vector field.}
 \label{fig:unsteadyNS_setting}
\end{figure}

We now consider the time-dependent case detailed in Equation~\eqref{eq:NS} with kinematic viscosity $\nu = 1 m^2 s^{-1}$ and $\left\lVert\mathbf{v}_{\mathrm{in}}\right\rVert = 10m s^{-1}$ in order to focus on a low Reynolds number ($\operatorname{Re} = 30$). The scenario parameters are therefore the inflow angle $\alpha_{\mathrm{in}}$
and, as always happens in time-dependent settings, the time variable $t$. We generate $22$ time-dependent optimal trajectories through the full-order OCP solver based on FEM implemented in \texttt{dolfin-adjoint}, solving the forward problem in Equation~\eqref{eq:NS} with the incremental Chorin-Temam projection method. Every trajectory considers a constant angle of attack randomly sampled in the interval $(-1.0, 1.0)$ radians and it is recorded in $10$ different time steps uniformly distributed in the interval $[0.05, 0.5]$ seconds, i.e. the final time is $T=0.5$ seconds and the time step is equal to $0.05$ seconds. The snapshots are then shuffled and divided into training set ($180$ snapshots) and test set ($40$ snapshots). Every optimal trajectory, which consists of $10$ different snapshots in time, is computed solving the nonlinear time-dependent full-order OCP requiring on average a computational time of $20$ minutes. An example of optimal control, velocity and pressure related to $\alpha_{\mathrm{in}} = 0.45$ radians and $t=0.5$ seconds is provided in Figure~\ref{fig:unsteadyNS_setting}: for visualization purposes, the scale considered for the control vectors is $10$ times bigger with respect to the figures in Section~\ref{subsec:steadyNS}.

\begin{figure}
    \centering
    \subfloat{\includegraphics[scale = 0.5]{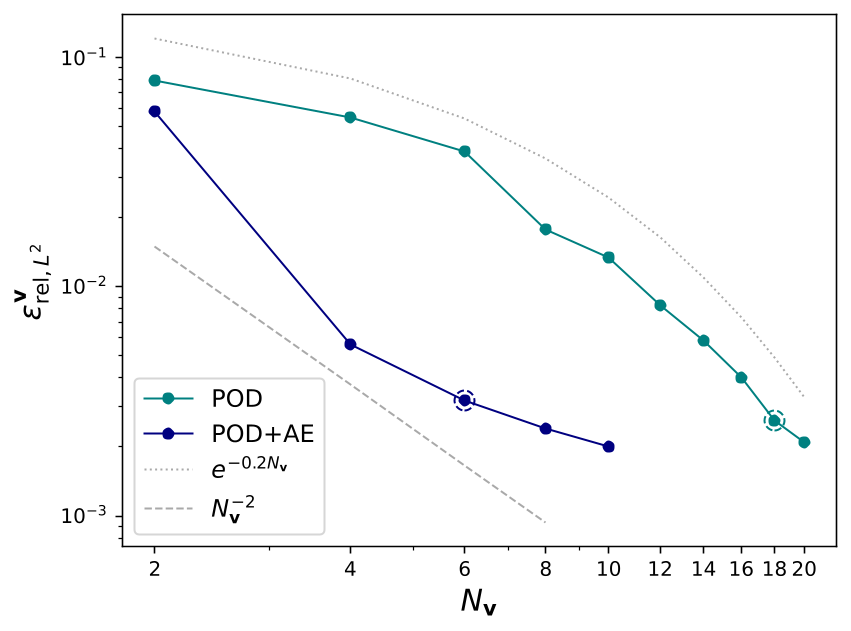}}
    \quad
    \subfloat{
    \includegraphics[scale = 0.5]{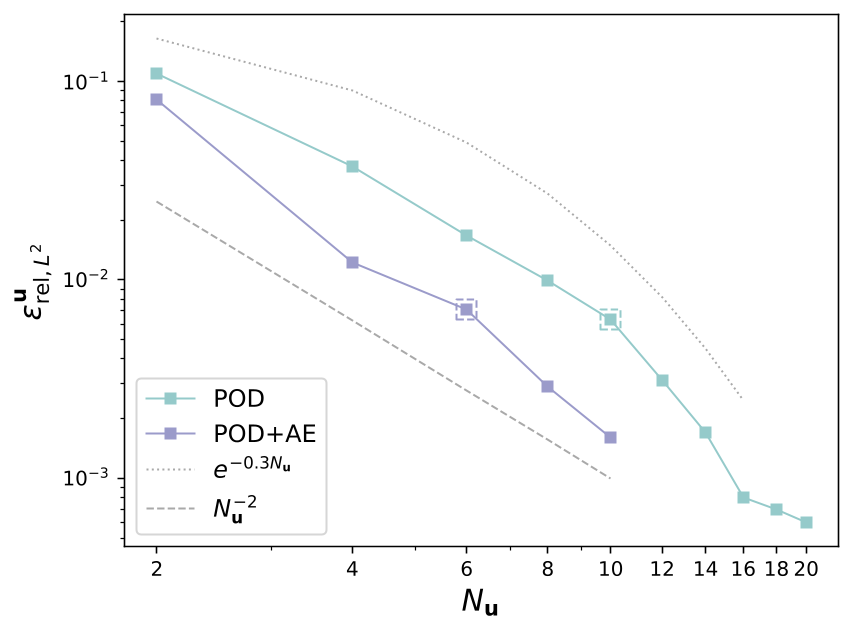}}

    \subfloat{\includegraphics[scale = 0.5]{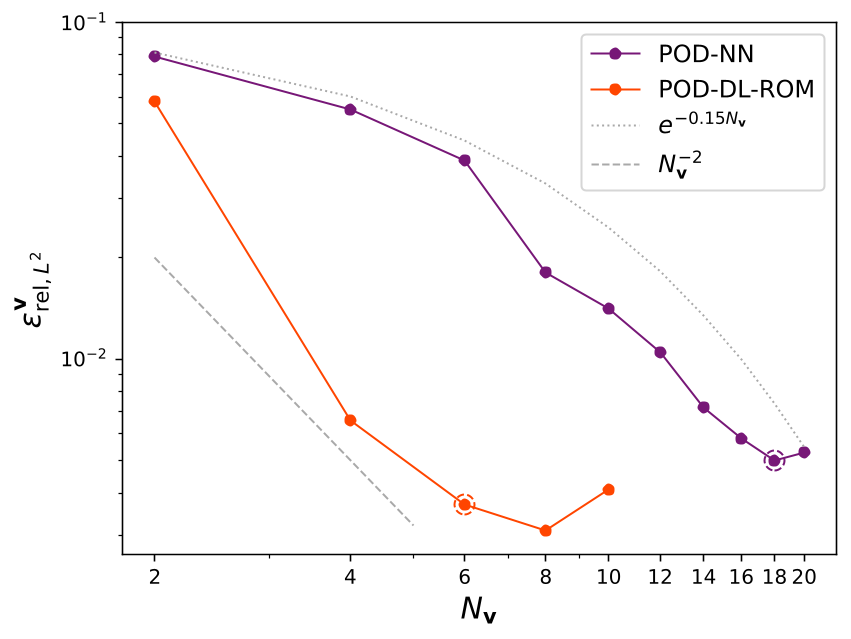}}
    \quad
    \subfloat{\includegraphics[scale = 0.5]{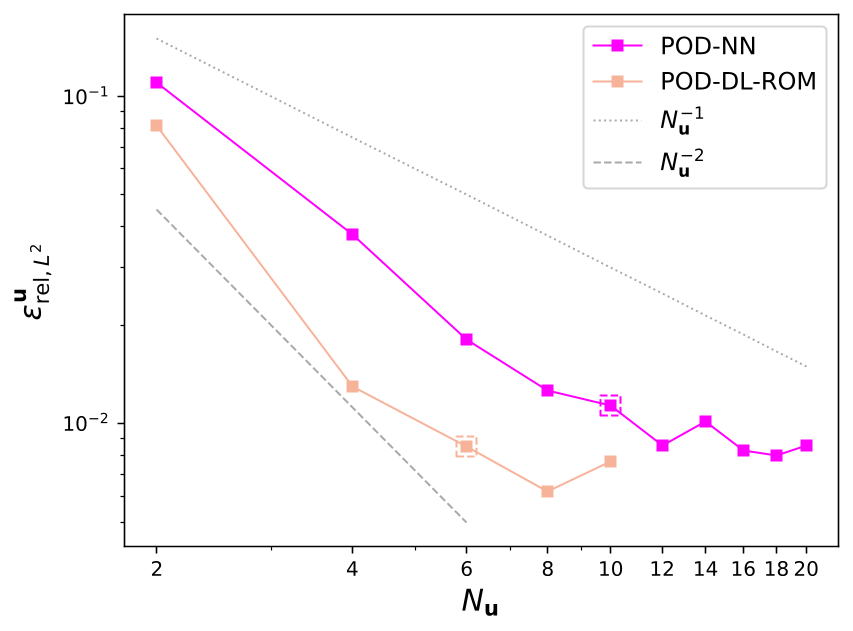}}
        
    \caption{\textit{Test 1.2}. Time-dependent flow control. Top: $L^2$ mean relative reconstruction error on $\mathbf{v}$ and $\mathbf{u}$ committed by POD and POD+AE for different latent dimensions $N_{\mathbf{v}}$ and $N_{\mathbf{u}}$. Bottom: $L^2$ mean relative prediction error on $\mathbf{v}$ and $\mathbf{u}$ committed by POD-NN and POD-DL-ROM for different latent dimensions $N_\mathbf{v}$ and $N_\mathbf{u}$. The latent dimensions selected in the numerical example are highlighted with dashed markers.}
   \label{fig:unsteadyNS_error_analysis1}
\end{figure}

In this setting, the $L^2$ mean relative test errors entailed when reconstructing state and control through POD shows an exponential decay, as visible in Figure~\ref{fig:unsteadyNS_error_analysis1}. However, due to a slow decay in correspondence of small latent dimensions, we need a higher number of POD modes with respect to the test case in Section~\ref{subsec:steadyNS} to achieve the same accuracy: specifically, with $N_{\mathbf{u}} = 10$ and $N_{\mathbf{v}} = 18$, the $L^2$ mean relative reconstruction errors on test data are equal to $0.26\%$ on the state and $0.63\%$ on the control.

\begin{figure}
\subfloat{\includegraphics[height = 0.3\linewidth]{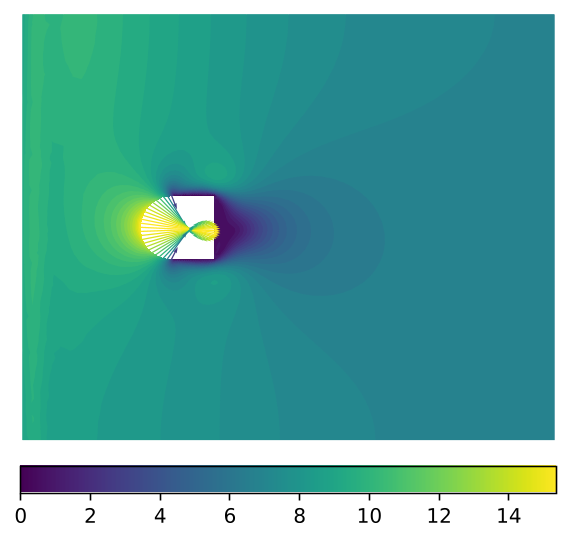}}
\subfloat{\includegraphics[height = 0.3\linewidth]{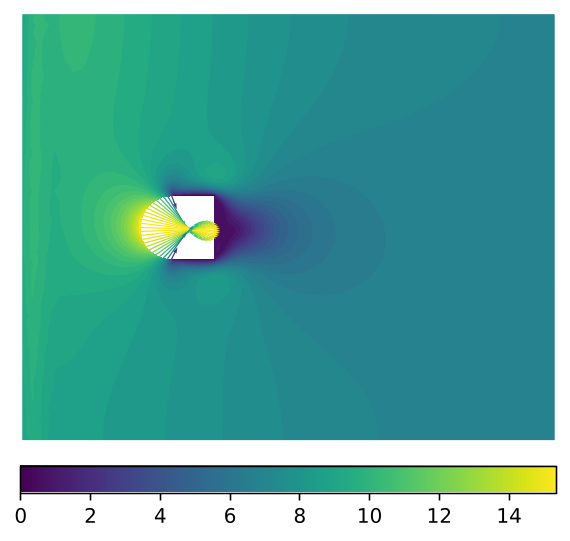}}
\subfloat{\includegraphics[height = 0.3\linewidth]{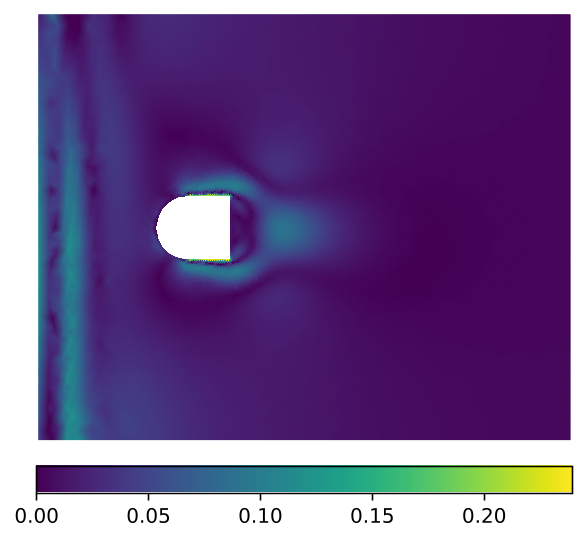}}

\hphantom{\subfloat{\includegraphics[height = 0.3\linewidth]{Images/unsteadyNS/Ground_truth_alpha_in=0.34_t=0.25.png}}}\subfloat{\includegraphics[height = 0.3\linewidth]{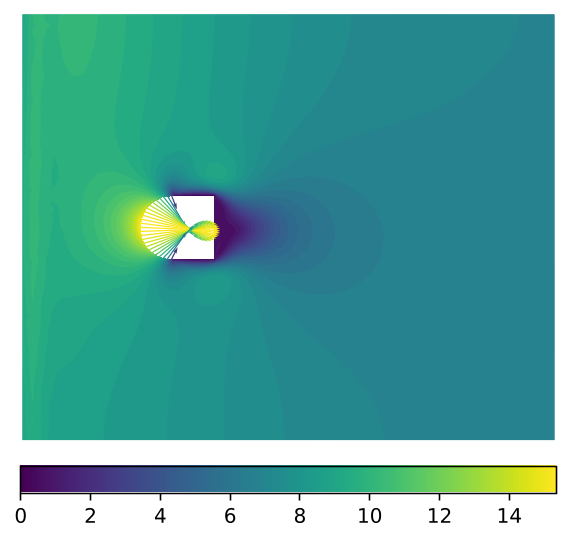}}
\subfloat{\includegraphics[height = 0.3\linewidth]{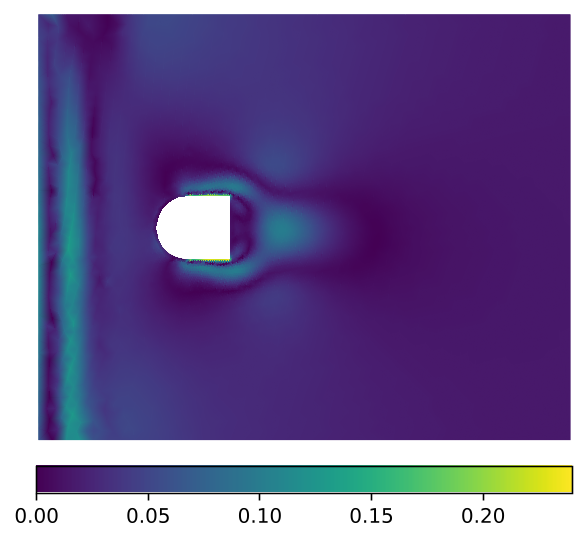}}

\caption{\textit{Test 1.2}. Time-dependent flow control. First row: high-fidelity optimal snapshot, POD+AE reconstruction and reconstruction error corresponding to the test scenario parameters $t = 0.25$ seconds and $\alpha_{\mathrm{in}} = 0.34$ radians. Second row: POD-DL-ROM prediction and prediction error corresponding to the test scenario parameters $t = 0.25$ seconds and $\alpha_{\mathrm{in}} = 0.34$ radians. The velocity on $\Omega$ is depicted through a scalar field with colours corresponding to its norm, while the control on $\Gamma_c$ is represented through a vector field.}
\label{fig:unsteadyNS_1}
\end{figure}

\begin{figure}
\centering
\subfloat{\includegraphics[height = 0.29\linewidth]{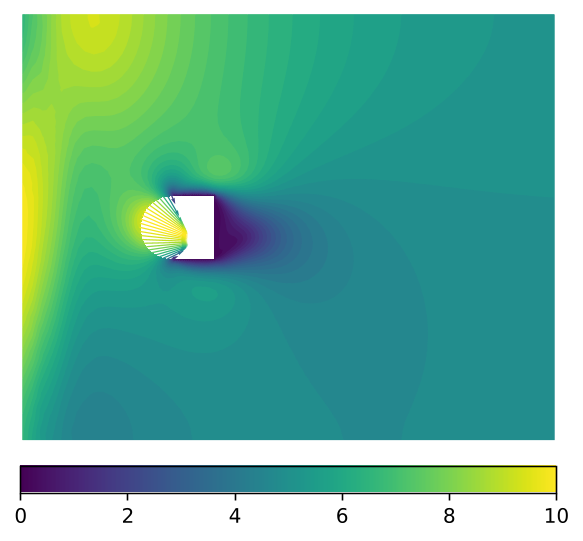}}
\subfloat{\includegraphics[height = 0.29\linewidth]{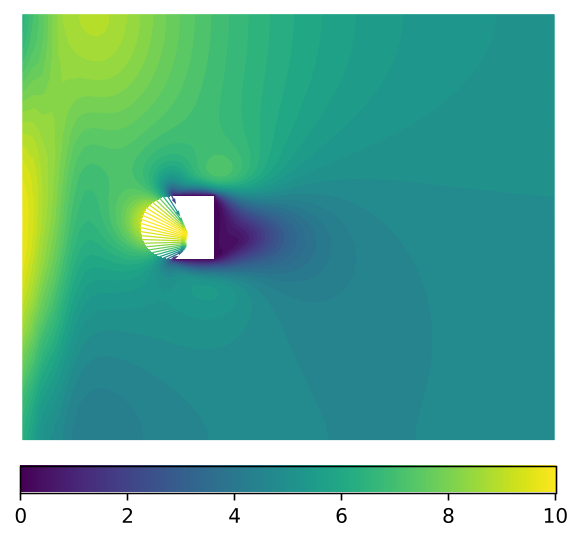}}
\subfloat{\includegraphics[height = 0.29\linewidth]{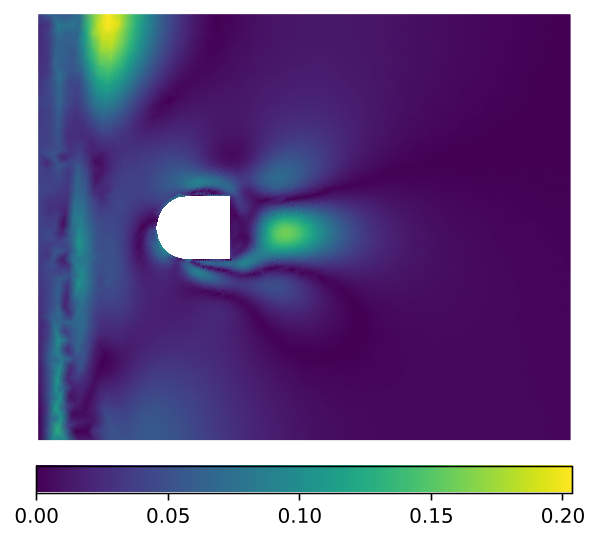}}

\hphantom{\subfloat{\includegraphics[height = 0.29\linewidth]{Images/unsteadyNS/Ground_truth_alpha_in=0.85_t=0.45.png}}}\subfloat{\includegraphics[height = 0.29\linewidth]{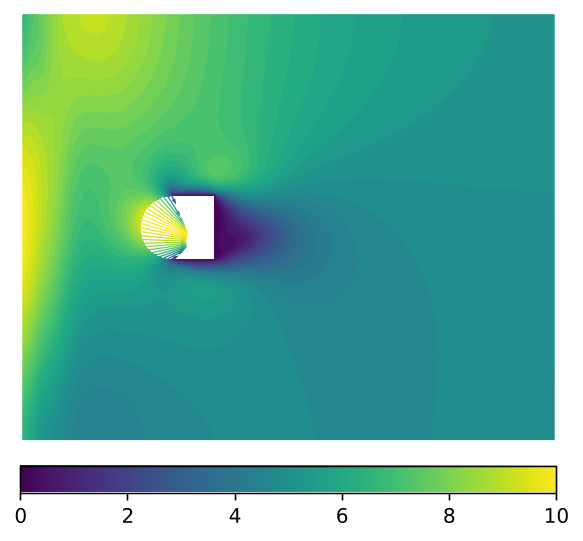}}
\subfloat{\includegraphics[height = 0.29\linewidth]{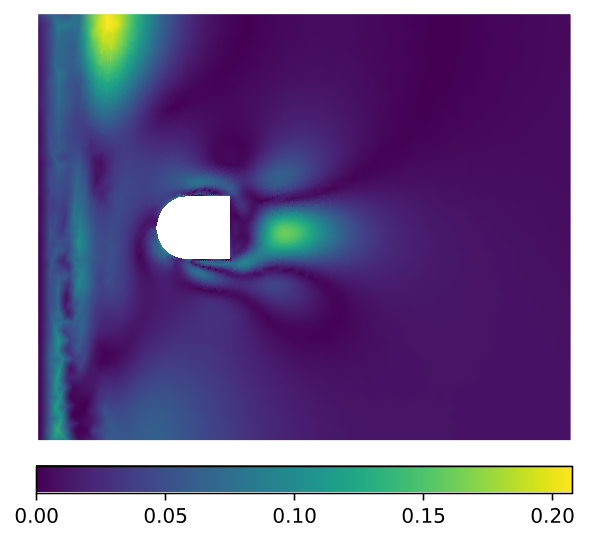}}

\caption{\textit{Test 1.2}. Time-dependent flow control. First row: high-fidelity optimal snapshot, POD+AE reconstruction and reconstruction error corresponding to the test scenario parameters $t = 0.45$ seconds and $\alpha_{\mathrm{in}} = 0.85$ radians. Second row: POD-DL-ROM prediction and prediction error corresponding to the test scenario parameters $t = 0.45$ seconds and $\alpha_{\mathrm{in}} = 0.85$ radians. The velocity on $\Omega$ is depicted through a scalar field with colours corresponding to its norm, while the control on $\Gamma_c$ is represented through a vector field.}
\label{fig:unsteadyNS_2}
\end{figure}

Nonlinear reduction strategies may be therefore considered to further reduce the problem dimensionality and allow for a lighter, more accurate and easily trainable neural network $\varphi$ with smaller output dimensions. As proposed by \citeinline{Fresca2022}{Fresca and Manzoni} and resumed in Section~\ref{sec:OCP-DL-ROM}, it is possible to combine an initial linear reduction step through POD and a nonlinear autoencoder to further compress the data. In particular, we consider $80$ and $140$ POD modes for the control and the state variables, respectively, which returns very low reconstruction errors on test data. Thanks to the nonlinear reduction provided by fully-connected autoencoders, the latent dimensions are then shrunk down to $N_{\mathbf{u}} = N_{\mathbf{v}} = 6$. The encoder responsible for the reduction of the state POD modes has $1$ hidden layer with $50$ neurons, while the corresponding decoder exploits $2$ hidden layers with $50$ neurons each. As far as the compression of control POD modes is concerned, the same number of hidden layers are considered, but the number of neurons per layer is decreased to $25$ thanks to the smaller input-output dimensions. Thanks to the joint reduction provided by POD and autoencoders, it is possible to achieve a very low-dimensional latent representation of the system variables with high-accuracy: indeed, the $L^2$ mean relative reconstruction errors on test data are $0.32 \%$ for the state and $0.70 \%$ for the control. The first rows of Figure~\ref{fig:unsteadyNS_1} and Figure~\ref{fig:unsteadyNS_2} display the ground truth and the reconstructed optimal pair related to two different test scenarios, confirming the high precision of the reduction procedure. The first row of Figure~\ref{fig:unsteadyNS_error_analysis1} shows, instead, the $L^2$ mean relative errors committed when reconstructing test data by POD and POD+AE, that is the combination of POD and autoencoders, for different latent dimensions, while keeping fixed the other hyperparameters. It is possible to assess that a nonlinear reduction technique allows to achieve more accurate results while considering smaller latent dimensions, entailing a lighter $\varphi$ architecture that is faster to train and to evaluate online.

\begin{figure}
    \centering
    \subfloat{\includegraphics[scale = 0.5]{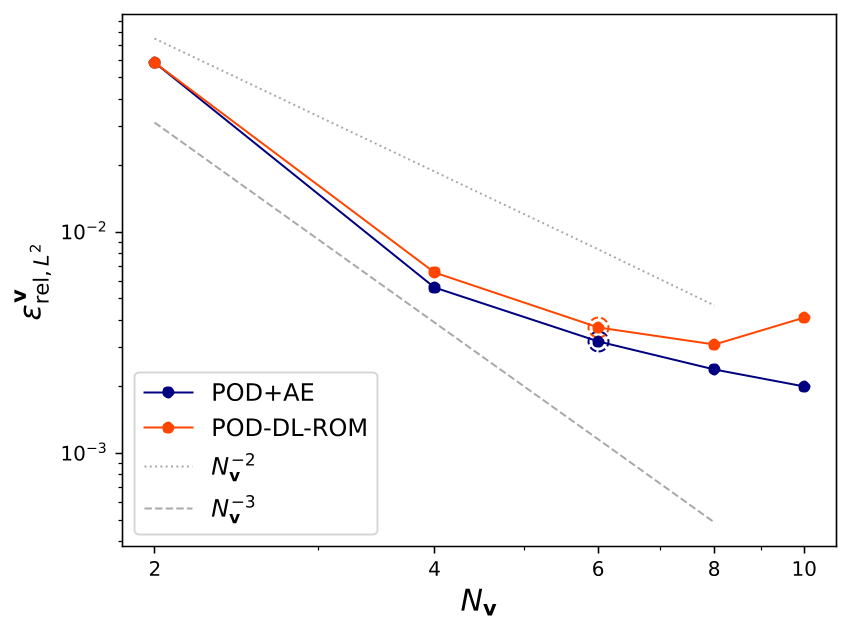}}
    \quad
    \subfloat{\includegraphics[scale = 0.5]{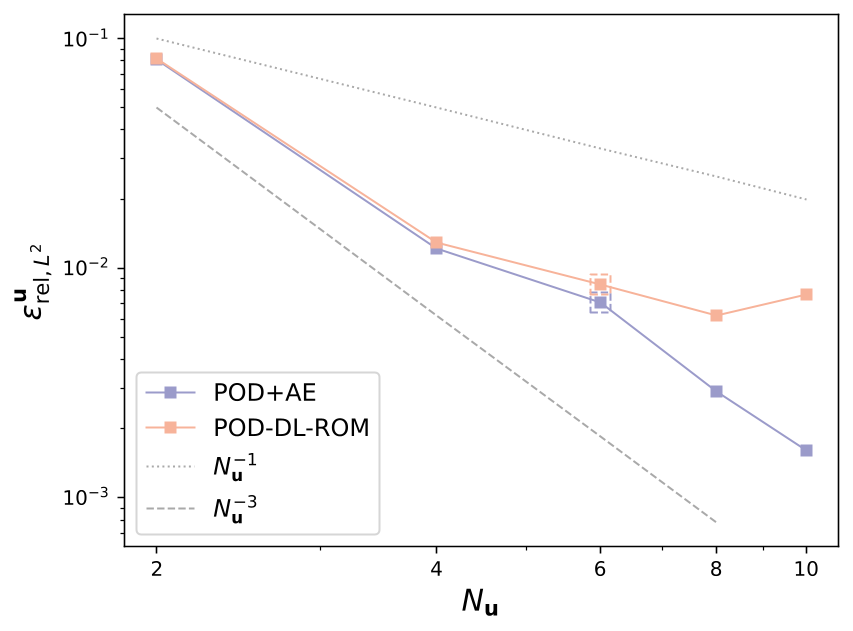}}
    
    \caption{\textit{Test 1.2}. Time-dependent flow control. Left: $L^2$ mean relative error on $\mathbf{v}$ committed by POD+AE (reconstruction error) and POD-DL-ROM (prediction error) for different latent dimensions $N_{\mathbf{v}}$. Right: $L^2$ mean relative error on $\mathbf{u}$ committed by POD+AE (reconstruction error) and POD-DL-ROM (prediction error) for different latent dimensions $N_{\mathbf{u}}$. The latent dimensions selected in
    the test case are highlighted with dashed markers.}
   \label{fig:unsteadyNS_error_analysis2}
\end{figure}

To rapidly compute the optimal state and control for different time instants and inflow velocities, we model the map from these parameters onto the reduced optimal pairs through a feed-forward neural network $\varphi$, thus considering a POD-DL-ROM approach. In particular, $\varphi$ consists of $2$ hidden layers with $150$ neurons each, and leaky Relu is exploited as activation function. After training the neural network with the L-BFGS optimization algorithm, the $L^2$ mean relative errors on test data are equal to $0.37\%$ for the optimal state and $0.85\%$ for the optimal control. The computational time required to train the autoencoders and the neural network $\varphi$ amounts to $4$ minutes and $11.52$ seconds. The second row of Figure~\ref{fig:unsteadyNS_error_analysis1} compares the $L^2$ mean relative errors committed when predicting test data through POD-NN and POD-DL-ROM for different latent dimensions: coherently with the corresponding reconstruction errors, a POD-DL-ROM allows to obtain more accurate results with smaller latent dimensions. Moreover, Figure~\ref{fig:unsteadyNS_error_analysis2} shows the reconstruction and prediction errors of POD-DL-ROM: while both show a superlinear decay with respect to the latent dimension, the prediction errors are slightly higher than the reconstruction ones due to the approximation of the parameter-to-solution map through the neural network $\varphi$. The second rows of Figure~\ref{fig:unsteadyNS_1} and Figure~\ref{fig:unsteadyNS_2} display, instead, the predictions obtained by POD-DL-ROM compared to the ground truth optimal solutions related to two different choices of test scenario parameters: both the control and the velocity are retrieved with high accuracy and no substantial differences can be assessed between ground truth data and POD-DL-ROM predictions. 

\begin{figure}
    \centering
    \subfloat{\includegraphics[height = 0.3\linewidth]{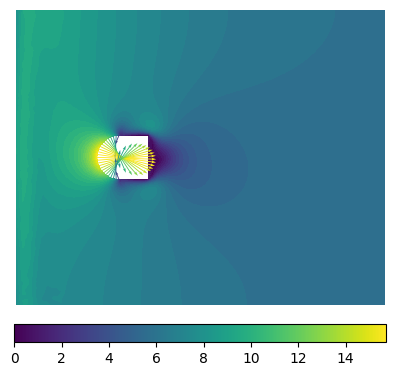}}
    \quad
    \subfloat{\includegraphics[height = 0.3\linewidth]{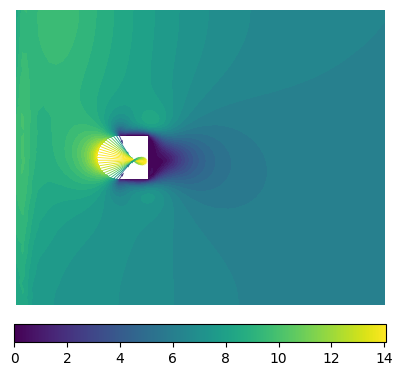}}
    \quad
    \subfloat{\includegraphics[height = 0.3\linewidth]{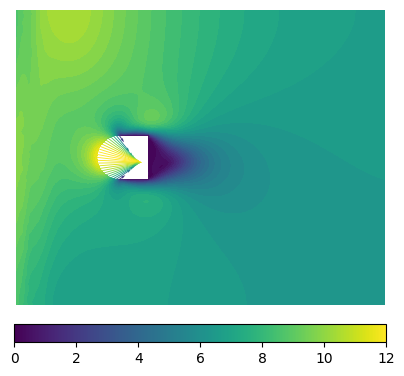}}

  \caption{\textit{Test 1.2}. Time-dependent flow control. Optimal state and control provided by POD-DL-ROM corresponding to the scenario parameters $\alpha_{\mathrm{in}} = 0.5$ radians and $t = 0.159, 0.276, 0.453$ seconds.}
\label{fig:unsteadyNS_results}
\end{figure}

The POD-DL-ROM architecture is then exploited to solve the time-dependent parametrized OCP for new and unseen time instants and angles of attack sampled in the intervals $[0.05, 0.5]$ seconds and $(-1.0, 1.0)$ radians, respectively. Figure~\ref{fig:unsteadyNS_results} shows the time-dependent optimal pairs computed via POD-DL-ROM when considering $\alpha_{\mathrm{in}} = 0.5$ radians at three different time instants. Note that, in order to find each of these solutions, we just need a forward pass of the neural network $\varphi$ and a decoding step through the autoencoders and POD that lasts at most $0.01$ seconds. Note also that, the forward pass $\varphi_D(\varphi(t^{\mathrm{new}}, \mus^{\mathrm{new}}))$ infers the time-dependent optimal pair -- i.e. the minimizer of the loss function computed in the interval $(0,T)$ with fixed final time $T$ -- evaluated at $t = t^{\mathrm{new}}$. 

\subsection{Active thermal cooling}
\label{subsec:cooling}

\begin{figure}
\centering
\begin{minipage}[c]{0.32\linewidth}
\centering
    \includegraphics[height = 0.8\linewidth]{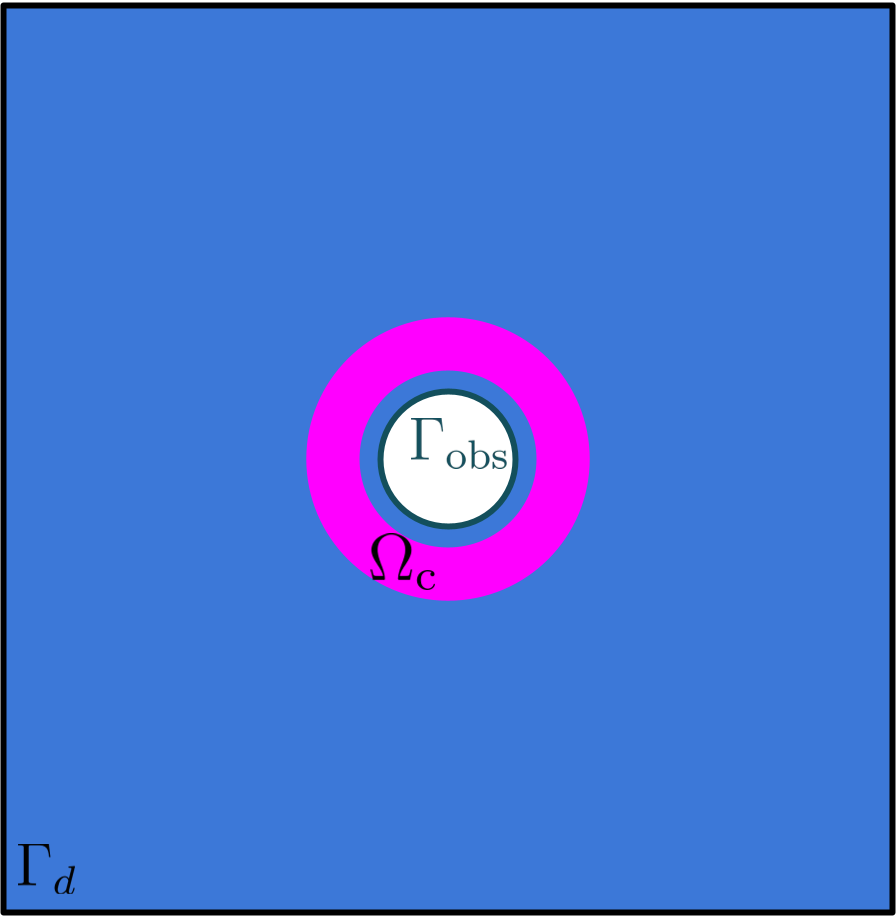}
\end{minipage}
\begin{minipage}[c]{0.66\linewidth}
    \subfloat{\includegraphics[height = 0.4\linewidth]{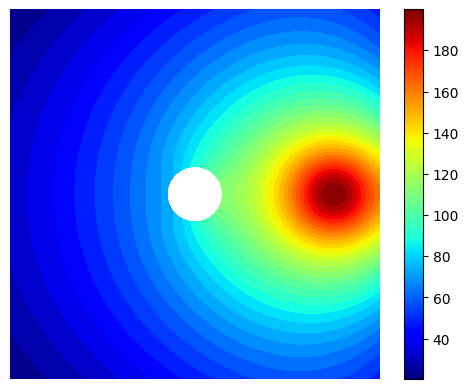}}
    \subfloat{\includegraphics[height = 0.4\linewidth]{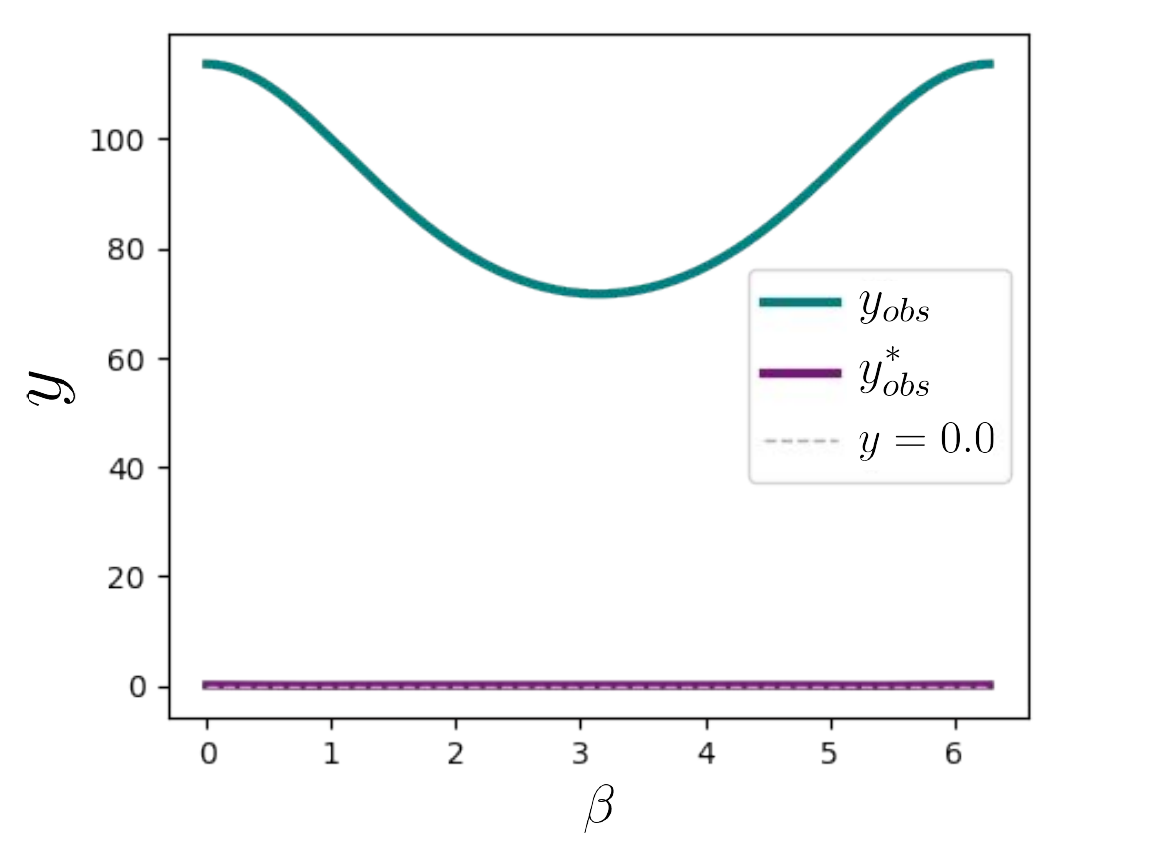}}

    \subfloat{\includegraphics[height = 0.4\linewidth]{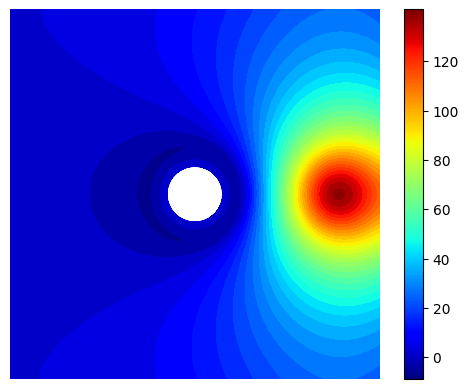}}
    \subfloat{\includegraphics[height = 0.4\linewidth]{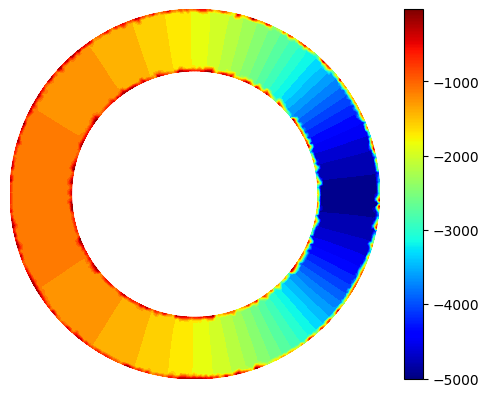}}
\end{minipage}
\caption{\textit{Test 2}. Active thermal cooling. Left: subdomains and boundaries considered in the active thermal cooling problem. In particular, $\Gamma_d$ in black is the external boundary, $\Gamma_{\mathrm{obs}}$ in green is the obstacle boundary and $\Omega_c$ in magenta is the annular region around the obstacle where the control source takes place. Center: temperature fields computed in absence of control action (top) and considering optimal control values (bottom) with thermal source located at $(0.75, 0.0)$. Top right: comparison of the temperature profiles at the obstacle boundary points $(0.15\cos(\beta),0.15\sin(\beta))$ for different angles $\beta$ in the controlled ($\mathit{y}_{obs}^*$) and uncontrolled ($y_{obs}$) setups. Bottom right: space-varying optimal control action considering the thermal source located at $(0.75, 0.0)$.}
\label{fig:cooling_setting}
\end{figure}

This section presents the numerical results obtained by applying the proposed real-time OCP solver to an active thermal cooling problem. We consider a square domain $\Omega = (-1,1) \times (-1,1)$ properly discretized thanks to \texttt{gmsh} utilities with the presence of a circular obstacle centered at $(0,0)$ with radius equal to $0.15$. The physical system under investigation is described by the steady heat equation, that is
\begin{equation}
    \begin{cases}
       - \nu \Delta y(\mathbf{x}) = s(\mathbf{x};\mus) + u(\mathbf{x})\mathbbm{1}_{\Omega_c}  \qquad &\mathrm{in} \ \Omega \\
       - \nu \nabla y(\mathbf{x}) \cdot \mathbf{n}(\mathbf{x}) = y(\mathbf{x})  \qquad &\mathrm{on} \ \Gamma_{\mathrm{d}} \\
        - \nu \nabla y(\mathbf{x}) \cdot \mathbf{n}(\mathbf{x}) = \gamma (y(\mathbf{x}) - y_{\mathrm{ext}})  \qquad &\mathrm{on} \ \Gamma_{\mathrm{obs}}
\end{cases}
\label{eq:cooling}
\end{equation}
where $\nu = 1 m^2 s^{-1}$ is the fixed material-specific thermal diffusivity and $\mathbf{n}$ is the unit vector normal to the domain boundary. The state variable $y$ is a scalar-quantity denoting the temperature (Kelvin degrees) at every domain location: since linear finite elements are considered, the vector $\yh$ collects the temperature values at these mesh nodes and in this case the dimension of the state space is $N_h^y = 51665$. The obstacle is heated by an external source $s$ ($K s^{-1}$) modeled through the Gaussian function
\[
s = 5000 \exp\{-40[(x_1-x_1^s)^2 + (x_2-x_2^s)^2]\}
\]
centered at the point $(x_1^s, x_2^s)$. To allow heat exchange between the obstacle and the surrounding environment, the Newton's law of cooling is taken into account as boundary condition on the obstacle boundary $\Gamma_{\mathrm{obs}}$, where $y_{\mathrm{ext}} = 125 K$ is the external temperature, that is proportional to the product between the intensity and the variance of $s$, and $\gamma = 1 m s^{-1}$ is the fixed heat transfer coefficient multiplied by the volumetric heat capacity. Moreover, a homogeneous Robin boundary condition is applied on the external boundary $\Gamma_{\mathrm{d}}$ in order to approximate an unbounded domain to first-order, as proposed by \citeinline{Sinigaglia2022-ROM}{Sinigaglia et al}. The problem setting is visible in the left panel of Figure~\ref{fig:cooling_setting}. 

In this context, our objective is to keep the temperature on the obstacle constant and equal, for simplicity, to the reference value $0 K$ while considering different locations of the external heat source $s$. The scenario parameters are therefore the coordinates of the source center: for convenience, the polar coordinates are taken into account, that is $\mus = [\vartheta_s, r_s]^{\top}$, where $x_1^s = r_s \cos(\vartheta_s)$ and $x_2^s = r_s \sin(\vartheta_s)$. To solve this parametric OCP, we consider an active control action, that may be interpreted as a cooler. In particular, $u$ is a space-varying source term taking values in $\Omega_c$, that is the annular region with inner and outer radius equal respectively to $0.2$ and $0.3$, as visible in the left panel of Figure~\ref{fig:cooling_setting}. The control variable is strongly high-dimensional since $N_h^u = 7721$ degrees of freedom are exploited to discretize the control region. In order to achieve the target temperature on $\Gamma_{\mathrm{obs}}$ while considering feasible and regular control actions, the cost function is defined as
\[
J(y(\mathbf{x}),u(\mathbf{x})) = \int_{\Gamma_{\mathrm{obs}}} \left|y(\mathbf{x})\right|^2 d\Gamma_{\mathrm{obs}} + \beta \int_{\Omega_c} \left|u(\mathbf{x})\right|^2 d\Omega_c + \beta_g \int_{\Omega_c} \left\lVert\nabla u(\mathbf{x})\right\rVert^2 d\Omega_c
\]
where the coefficients $\beta$ and $\beta_g$ are set equal to $1\mathrm{e}{-8}$. To better understand the problem at hand, Figure~\ref{fig:cooling_setting} shows a comparison between the uncontrolled and the controlled temperature fields for a particular scenario: it is possible to assess that an optimal space-varying control strategy allows to cool the obstacle and keep its temperature equal to the reference value.

As far as the data generation process is concerned, $100$ optimal snapshots are computed through the OCP high-fidelity solver based on FEM provided by \texttt{dolfin-adjoint} considering random source coordinates $\vartheta_s, r_s$ uniformly sampled in the intervals $(-\frac{\pi}{2}, \frac{\pi}{2})$ radians and $(0.4, 0.9)$. The full-order solver requires, on average, $4$ minutes and $22$ seconds to retrieve one optimal snapshot while exploiting L-BFGS as optimization algorithm and a tolerance equal to $1\mathrm{e}{-10}$. The data are then divided into training set ($80$ snapshots) and test set ($20$ snapshots) for model evaluation purposes.

\begin{figure}
\centering
\begin{minipage}[c]{0.33\linewidth}
\subfloat{\includegraphics[height = 0.75\linewidth]{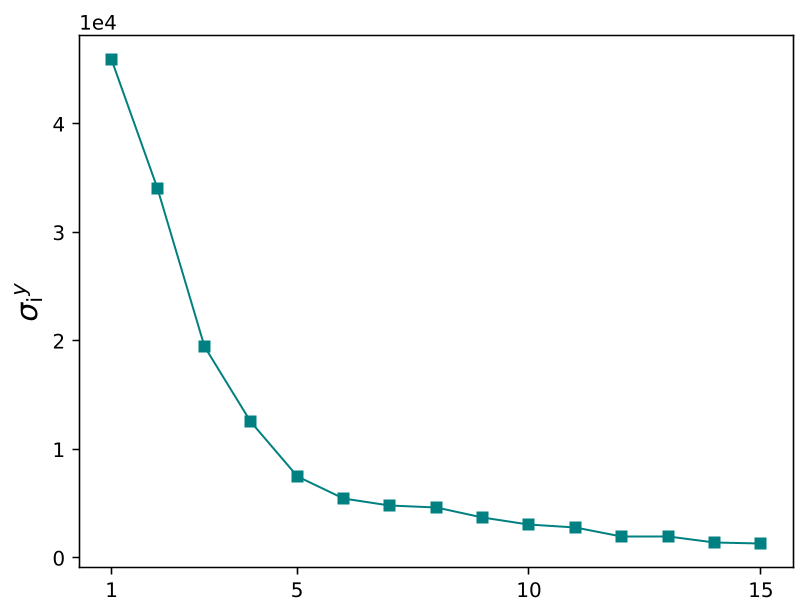}}

\subfloat{\includegraphics[height = 0.75\linewidth]{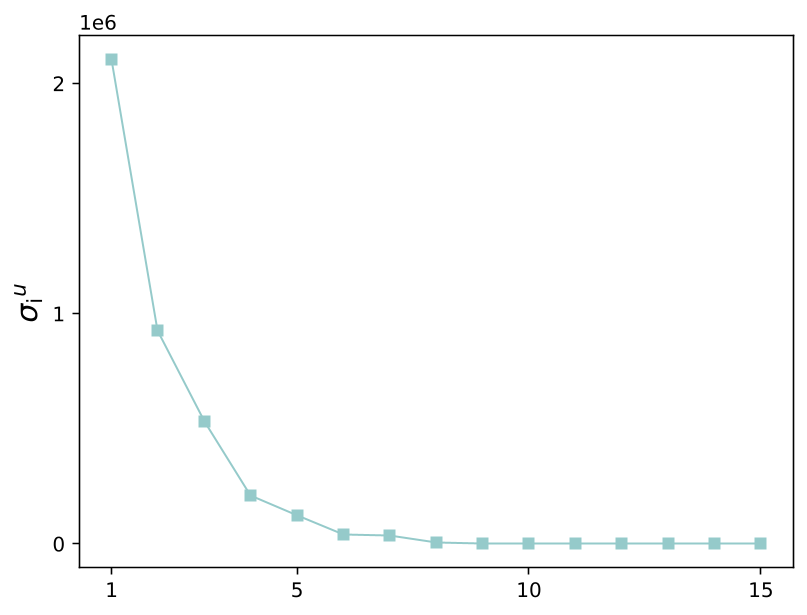}}
\end{minipage}
\begin{minipage}[c]{0.33\linewidth}
\subfloat{\includegraphics[height = 0.75\linewidth]{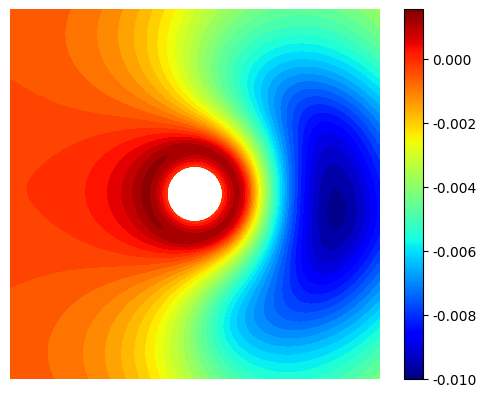}}

\subfloat{\includegraphics[height = 0.75\linewidth]{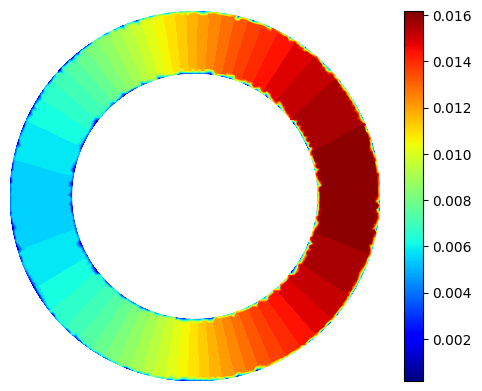}}
\end{minipage}
\begin{minipage}[c]{0.33\linewidth}
\subfloat{\includegraphics[height = 0.75\linewidth]{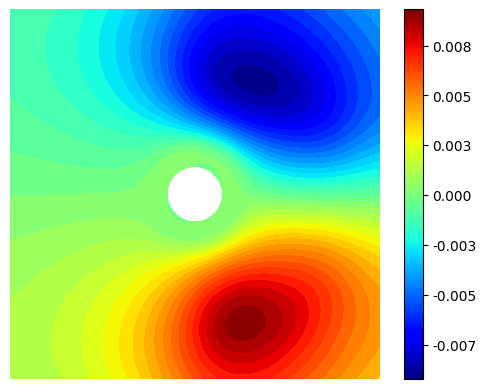}}

\subfloat{\includegraphics[height = 0.75\linewidth]{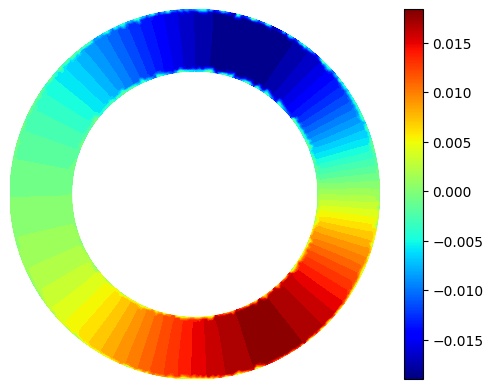}}
\end{minipage}

  \caption{\textit{Test 2}. Active thermal cooling. Left: decay of the singular values computed through SVD applied to the state (top) and control (bottom) snapshots matrices. Center: POD mode associated with the highest singular value of the state (top) and control (bottom). Right: POD mode associated with the second highest singular value of the state (top) and control (bottom).}
\label{fig:cooling_pod}
\end{figure}

In this context, dimensionality reduction is strongly recommended due to the huge dimensionality of state and control variables. In particular, a linear reduction scheme based on POD is enough to reduce the optimal control: indeed, considering $N_u = 7$, the $L^2$ mean relative reconstruction error obtained on test data results equal to $0.17\%$. The second row of Figure~\ref{fig:cooling_pod} displays the singular value decay and the two most-energetic POD modes -- i.e. the ones associated with the two largest singular values -- that capture the most of the control snapshots variance.

\begin{figure}
\centering
\subfloat{\includegraphics[height = 0.25\linewidth]{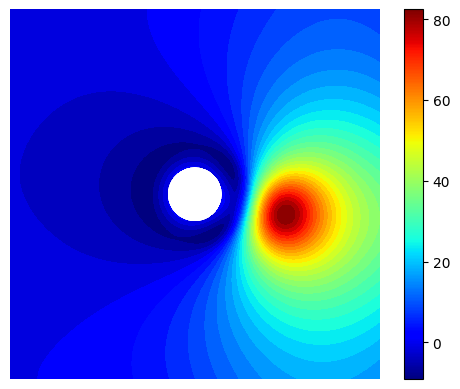}}
\subfloat{\includegraphics[height = 0.25\linewidth]{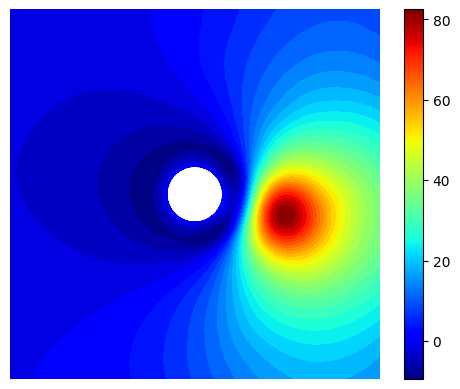}}
\subfloat{\includegraphics[height = 0.25\linewidth]{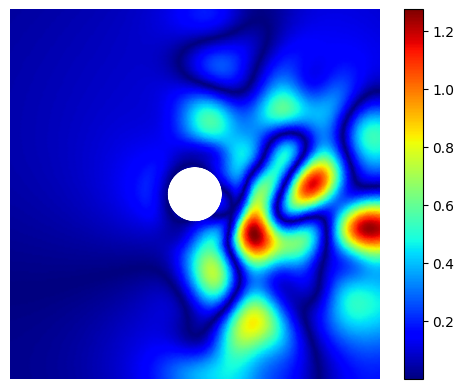}}

\hphantom{\subfloat{\includegraphics[height = 0.25\linewidth]{Images/cooling/Ground_truth_state_theta=-0.22_r=0.44.png}}}\subfloat{\includegraphics[height = 0.25\linewidth]{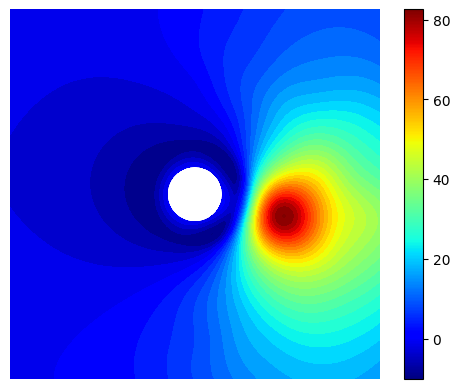}}
\subfloat{\includegraphics[height = 0.25\linewidth]{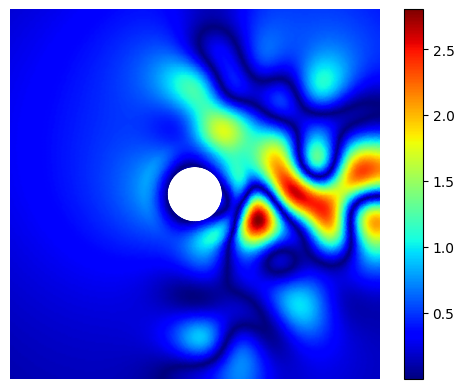}}

\caption{\textit{Test 2}. Active thermal cooling. First row: high-fidelity optimal snapshot, POD+AE reconstruction and reconstruction error corresponding to the test scenario parameters $\vartheta_s=-0.22$ radians and $r_s=0.44$. Second row: POD-DL-ROM prediction and prediction error
corresponding to the test scenario parameters $\vartheta_s=-0.22$ radians and $r_s=0.44$.}
\label{fig:cooling_state}
\end{figure}

\begin{figure}
\centering
\subfloat{\includegraphics[height = 0.25\linewidth]{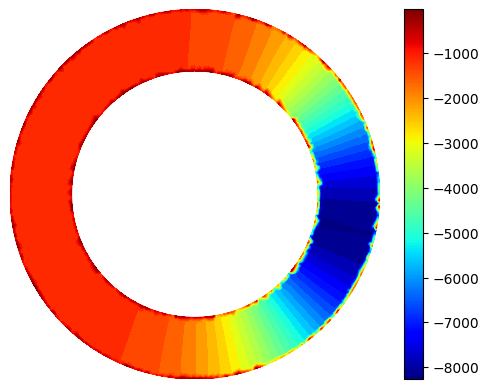}}
\subfloat{\includegraphics[height = 0.25\linewidth]{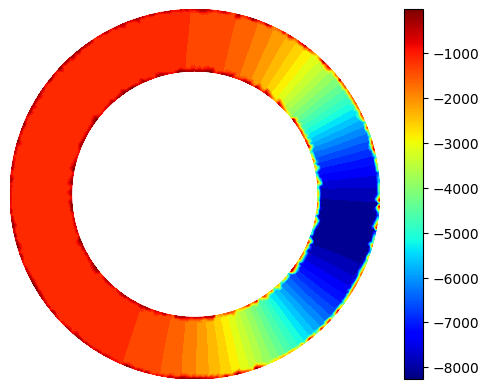}}
\subfloat{\includegraphics[height = 0.25\linewidth]{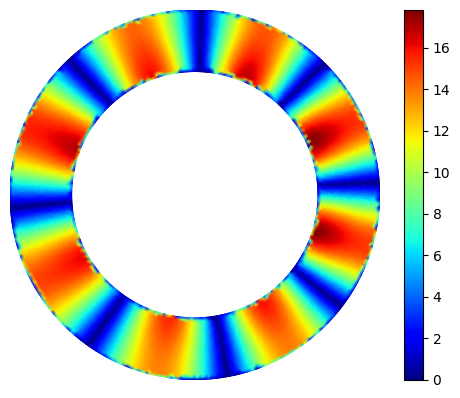}}

\hphantom{\subfloat{\includegraphics[height = 0.25\linewidth]{Images/cooling/Ground_truth_control_theta=-0.22_r=0.44.png}}}\subfloat{\includegraphics[height = 0.25\linewidth]{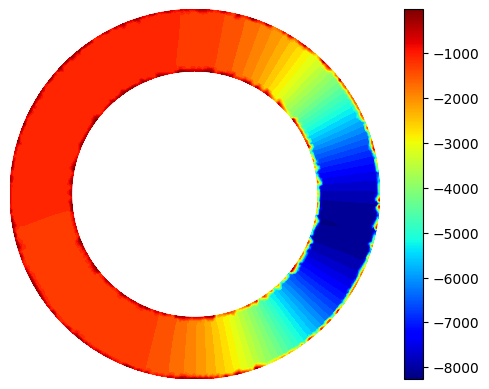}} 
\subfloat{\includegraphics[height = 0.25\linewidth]{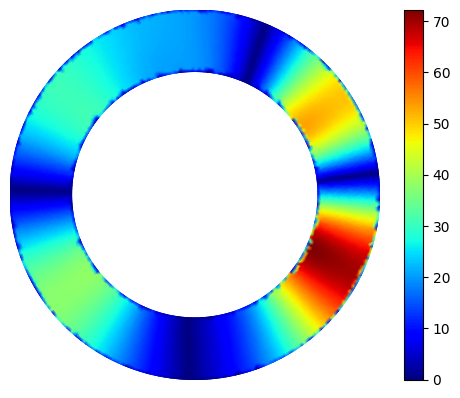}}

\caption{\textit{Test 2}. Active thermal cooling. First row: high-fidelity optimal snapshot, POD reconstruction and reconstruction error corresponding to the test scenario parameters $\vartheta_s=-0.22$ radians and $r_s=0.44$. Second row: POD-DL-ROM prediction and prediction error
corresponding to the test scenario parameters $\vartheta_s=-0.22$ radians and $r_s=0.44$.}
\label{fig:cooling_control}
\end{figure}

While POD is accurate to project the control snapshots in a low-dimensional latent space, the same is not true for state snapshots compression. Indeed, due to a slower singular values decay, as visible in Figure~\ref{fig:cooling_pod}, more than $37$ POD modes would be necessary in order to achieve acceptable reconstruction errors. Therefore, following the POD-DL-ROM strategy, $72$ POD modes are initially computed and then compressed thanks to a nonlinear autoencoder into a six-dimensional latent representation, that is $N_y = 6$. The encoder consists of $2$ layers with, respectively, $70$ and $30$ neurons, while $3$ layers having $30$, $50$ and $70$ neurons are considered in the decoder architecture, with leaky Relu as activation function. After training the autoencoder with the L-BFGS optimization algorithm, the $L^2$ mean relative reconstruction error committed on test data is equal to $2\%$. The first row in Figure~\ref{fig:cooling_state} displays the ground truth and the reconstructed temperature field corresponding to the test scenario parameters $\vartheta_s = -0.22$ radians and $r_s = 0.44$. The first row of Figure~\ref{fig:cooling_control} shows, instead, the ground truth and the corresponding optimal control compressed and reconstructed through POD related to the same choice of scenario parameters.

As detailed in Section~\ref{sec:OCP-DL-ROM}, the last ingredient required to obtain a non-intrusive real-time optimal solver for this parametric OCP is the parameter-to-solution map $\varphi$. Specifically, $\varphi$ is modeled through a feed-forward neural network having $2$ layers with $50$ neurons each and exploiting leaky Relu as activation function. The input layer of $\varphi$ is augmented considering, in addition to the scenario parameters $\vartheta_s$ and $r_s$, meaningful quantities such as $r_s\cos\vartheta_s$ and $r_s\sin\vartheta_s$. POD-DL-ROM is able to achieve a $L^2$ mean relative error on test data equal to $3.83\%$ for the optimal state and $1.10\%$ for the optimal control. The overall training of $\varphi$ and the autoencoder lasts $4$ minutes and $41$ seconds. Instead, the evaluation time required to predict the optimal pair related to a new scenario of interest is equal to, on average, $0.007$ seconds, thus providing a solver $40000$ times faster with respect to the full-order one. The bottom rows of Figure~\ref{fig:cooling_state} and Figure~\ref{fig:cooling_control} display a comparison between the ground truth optimal solutions $\mathbf{y}_h, \mathbf{u}_h$ and the POD-DL-ROM predictions obtained by
\[
\begin{bmatrix} \tilde{\mathbf{y}}_h \\ \tilde{\mathbf{u}}_h \end{bmatrix} = \varphi_D(\varphi(\vartheta_s, r_s))
\]
related to two different choices of scenario parameters in the test set, where $\varphi_D$ resumed the decoding actions applied separately to state and control.

\begin{figure}
\centering
\begin{minipage}[c]{0.33\linewidth}
     \subfloat{\includegraphics[height = 0.75\linewidth]{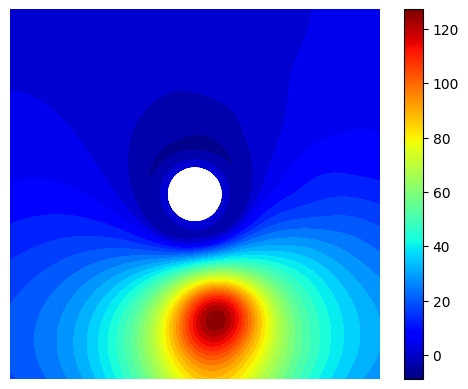}}

     \subfloat{\includegraphics[height = 0.75\linewidth]{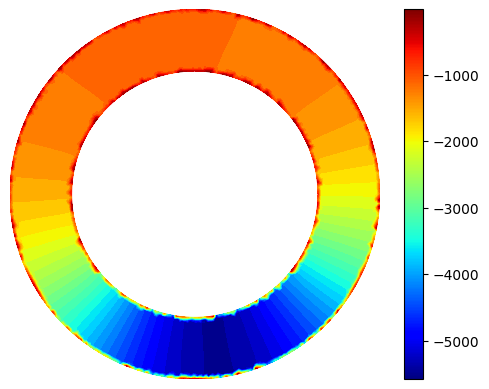}}
\end{minipage}
\begin{minipage}[c]{0.33\linewidth}
    \subfloat{\includegraphics[height = 0.75\linewidth]{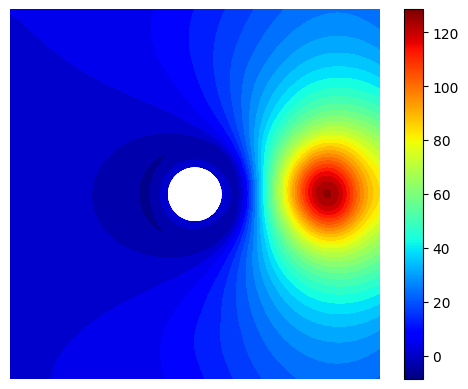}}
    
    \subfloat{\includegraphics[height = 0.75\linewidth]{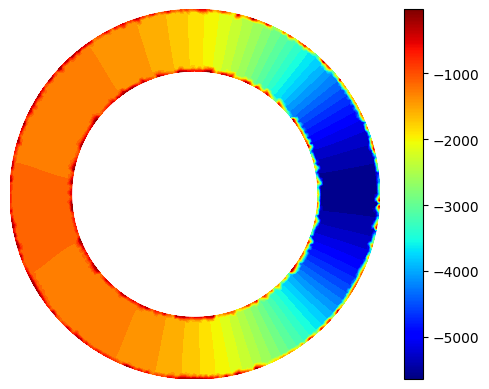}}
\end{minipage}
\begin{minipage}[c]{0.33\linewidth}
    \subfloat{\includegraphics[height = 0.75\linewidth]{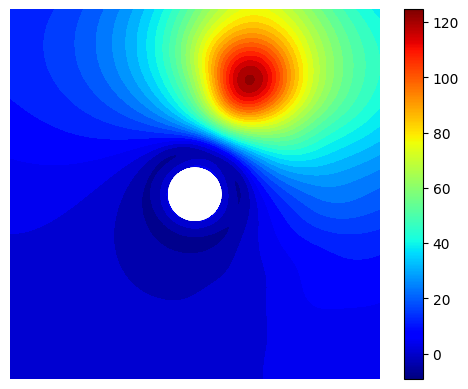}}
    
   \subfloat{\includegraphics[height = 0.75\linewidth]{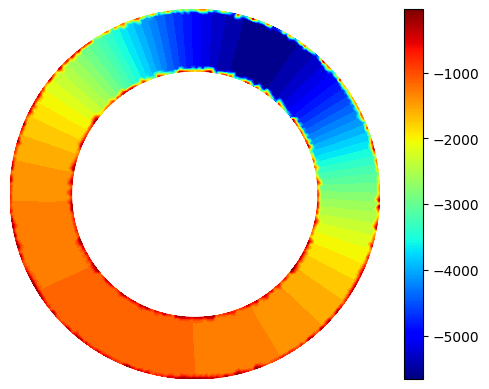}}
\end{minipage}

  \caption{\textit{Test 2}. Active thermal cooling. Optimal state and control provided by POD-DL-ROM corresponding to the scenario parameters $r_s = 0.65$ and $\vartheta_s = -1.42, 0.03, 1.13$ radians.}
\label{fig:cooling_test}
\end{figure}

\begin{figure}
    \centering
    \subfloat{\includegraphics[height = 0.34\linewidth]{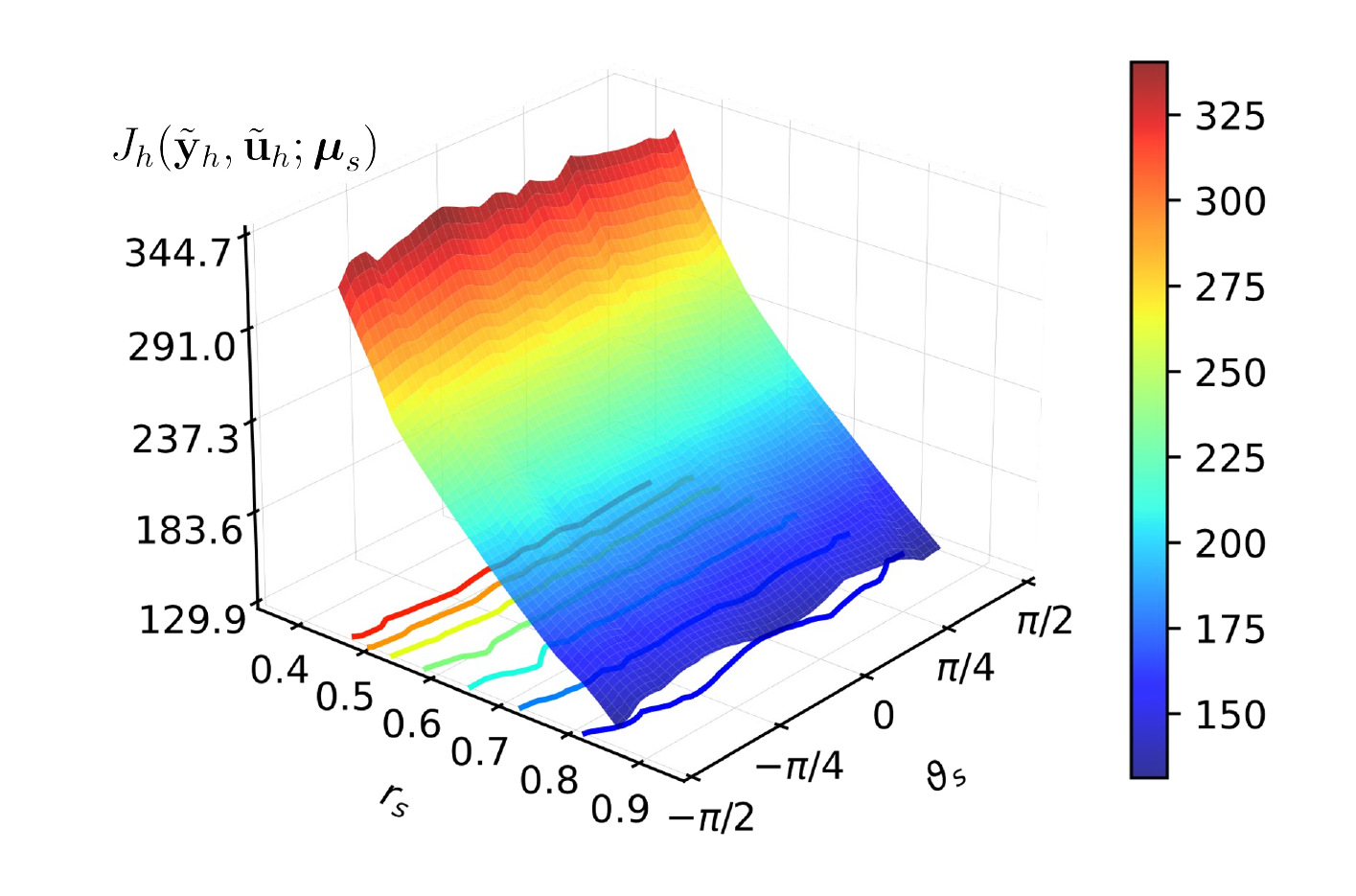}}
    \subfloat{\includegraphics[height = 0.34\linewidth]{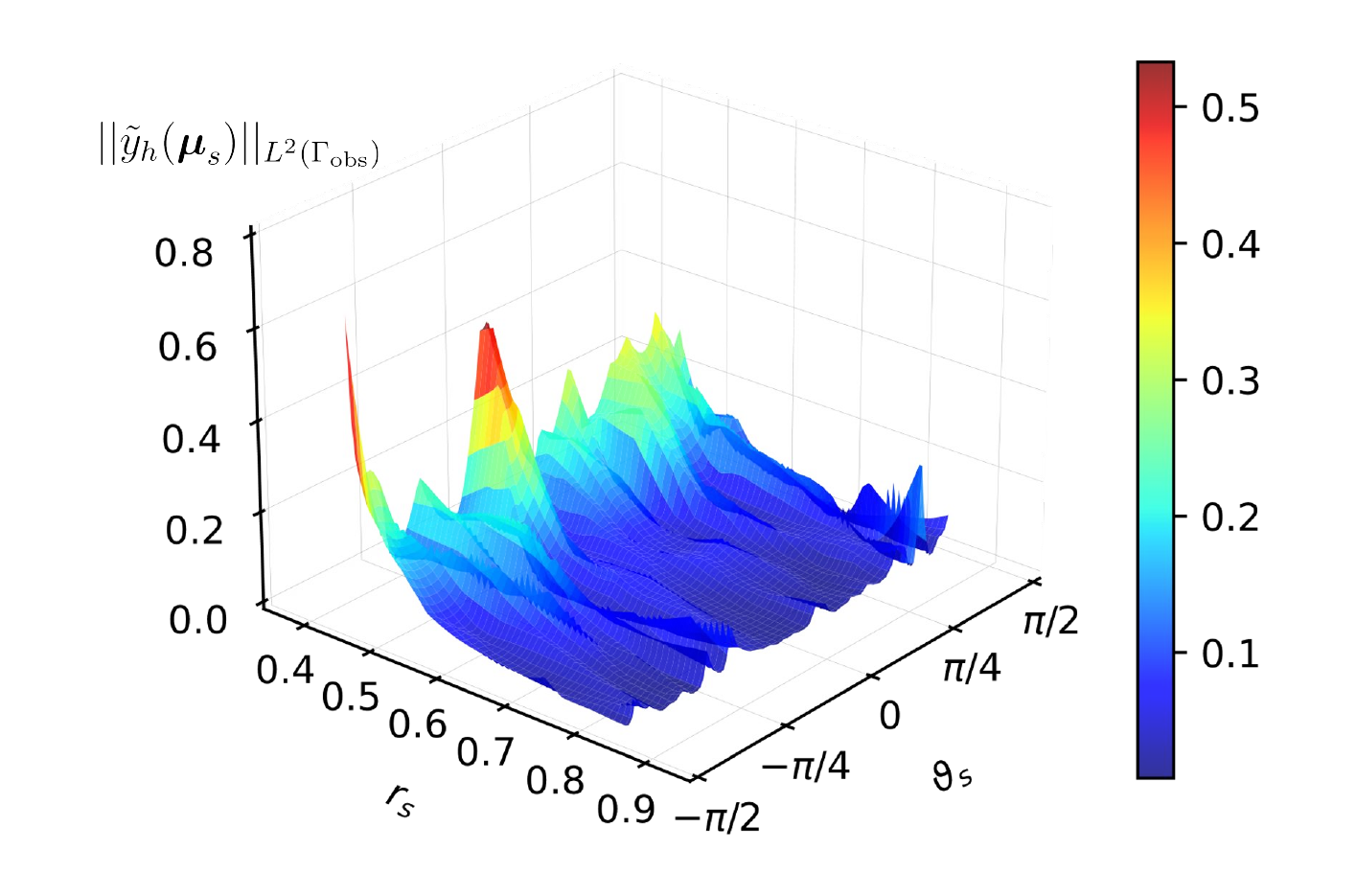}}

  \caption{\textit{Test 2}. Active thermal cooling. Left: optimal loss function values computed online considering $10000$ different scenario parameters in the parameter space. Contour lines of the loss function surface are available on the $\vartheta_s$-$r_s$ plane. Right: $L^2$ norm of the obstacle temperatures computed online considering $10000$ different scenario parameters in the parameter space.}
\label{fig:J}
\end{figure}

The POD-DL-ROM architecture is then tested online with respect to new scenario parameters unseen during training. Figure~\ref{fig:cooling_test} displays the results obtained for $r_s = 0.65$ and three different values of $\vartheta_s$. On average, the temperature on the obstacle is equal to $0.0011 K$, thus very close to the target value ($0K$). The left panel of Figure~\ref{fig:J} shows, instead, the optimal values of the loss function $J_h(\tildeyh, \tildeuh; \mus)$ evaluated at optimal states and controls computed for several scenarios in a uniform $100 \times 100$ grid covering the parameter space. Note that, since $10000$ optimal pairs are required in order to perform this parametric analysis, a full-order solver would need approximately $30$ days of computation. Instead, thanks to our fast-evaluable architecture, all the optimal solutions are computed in $2.14$ seconds. From the loss function behaviour in the parameter space, it is clear that the closer the heat source to the obstacle -- i.e. the smaller $r_s$ -- the higher the loss due to more energetic control terms required. The same analysis is performed looking at the optimal $L^2$ norm of the temperature on the obstacle boundary, that is \[||\tilde{y}_h(\mus)||_{L^2(\Gamma_{\mathrm{obs}})} = \sqrt{\int_{\Gamma_{\mathrm{obs}}} \tilde{y}_h^2(\mathbf{x};\mus) d\Gamma_{\mathrm{obs}}}\] where $\tilde{y}_h$ is the function associated with the FEM coefficients $\tildeyh$. As visible in the right panel of Figure~\ref{fig:J}, the control strategies predicted by POD-DL-ROM allow to attain obstacle temperatures close to the target in the whole parameter space and, thus, to properly cool the obstacle. In particular, since the loss function is a sum of terms being minimized, slightly higher obstacle temperatures are naturally obtained corresponding to more energetic control strategies.
 
\section{Conclusions}\label{sec:conclusions}

In this work we propose an efficient and reliable non-intrusive data-driven method to solve parametrized OCPs. Indeed, in all applications we considered, the mismatch between the predicted optimal solutions and the corresponding ground truth is always less than $4\%$ when taking into account new unseen scenarios. Moreover, the online stage is really fast -- it lasts at most $0.01$ seconds -- regardless of the problem complexity and dimension. Along with speed and accuracy, another key feature of the presented approach is its extreme versatility: in fact, as demonstrated throughout the numerical results in Section~\ref{sec:test}, it is possible to handle a broad spectrum of parametrized OCPs, ranging from scalar to vector, from linear to nonlinear and from steady to time-dependent problems.

An additional key point of the proposed strategy is the dimensionality reduction step that, differently from several control strategies available in the literature, allows us to easily deal with high-dimensional state variables and boundary or distributed controls. Furthermore, different non-intrusive reduced order models, including POD-NNs, DL-ROMs, and POD-DL-ROMs, can be taken into account in order to exploit the most effective reduction strategy for each OCP at hand, extending the current state-of-the-art on non-intrusive ROMs to control problems.

Several extensions of the proposed framework can be considered, even if they are left to future works. For instance, instead of dealing with high-dimensional states simulated through FEM solvers, it is in principle possible to exploit sensor data or images representing the underlying dynamics. A further, interesting enhancement of the proposed framework is to consider data-driven or physics-informed surrogate models approximating the forward problem to speed up the data generation process, which can be computationally intensive, although necessary only once during the offline phase. Moreover, uncertain parameters can be incorporated into the parameterized PDE to enable real-time robust control strategies. Feedback signals may also be considered if additional measurements are available online, as well as attention mechanisms may be integrated into the parameter-to-solution map in order to achieve better results when extrapolating in time.


\section*{Acknowledgments}
AM acknowledges the Project “Reduced Order Modeling and Deep Learning for the real-time approximation of PDEs (DREAM)” (Starting Grant No. FIS00003154), funded by the Italian Science Fund (FIS) - Ministero dell'Università e della Ricerca and the project FAIR (Future Artificial Intelligence Research), funded by the NextGenerationEU program within the PNRR-PE-AI scheme (M4C2, Investment 1.3, Line on Artificial Intelligence). MT and AM are members of the Gruppo Nazionale Calcolo Scientifico-Istituto Nazionale di Alta Matematica (GNCS-INdAM).

\bibliographystyle{abbrv}
\bibliography{references}  

\end{document}